\documentclass[11pt]{article}
\usepackage{amssymb}
\usepackage{mathrsfs}
\usepackage{amsfonts}
\usepackage{tikz}
\usepackage{titlesec}
\usepackage[titletoc]{appendix}
\usepackage{ifpdf}

\usepackage{amsmath,amsthm,amsfonts,amssymb,amscd}

\allowdisplaybreaks[4]
\newtheorem {Lemma}{Lemma}[section]
\newtheorem{Corollary}[Lemma]{Corollary}
\newtheorem {Theorem} [Lemma]{Theorem}

\newenvironment {Proof}{\noindent {\bf Proof.}}{\hfill\ensuremath{\square}}
\newcommand*{\QEDB}{\hfill\ensuremath{\square}}

\begin{document}

\title{Some new sufficient conditions for  $2p$-Hamilton-biconnectedness of graphs \thanks{This work is supported by  the Joint NSFC-ISF Research Program (jointly funded by the National Natural Science Foundation of China and the Israel Science Foundation (No. 11561141001)),  the National Natural Science Foundation of China (No.11531001 ).
}}

\author{ Ming-Zhu Chen, Xiao-Dong Zhang\footnote{Corresponding author. E-mail: xiaodong@sjtu.edu.cn}\\
School of Mathematics Science, MOE-LSC, SHL-MAC\\
Shanghai Jiao Tong University,
Shanghai 200240, P. R. China}

\date{}
\maketitle

\begin{abstract}

A   balanced bipartite  graph  $G$ is said to be \emph{$2p$-Hamilton-biconnected} if for any balanced subset $W$ of size $2p$ of $V(G)$, the subgraph induced by $V(G)\backslash W$ is Hamilton-biconnected. In this paper, we prove that
``
Let $p\geq0$ and $G$ be a  balanced bipartite graph of order $2n$ with minimum degree $\delta(G)\geq k$, where $n\geq 2k-p+2$ and  $k\geq p$. If
 the number of edges $
  e(G)>n(n-k+p-1)+(k+2)(k-p+1),
$ then $G$ is $2p$-Hamilton-biconnected except some exceptions.'' Furthermore, this result is used to present two new spectral conditions for a graph to $2p$-Hamilton-biconnected.
Moreover, the similar results are also presented for nearly balanced bipartite graphs.

\end{abstract}

\section{Introduction}
 Let $G$ be an   undirected  simple  graph with vertex set
$V(G)=\{v_1,\dots,v_n\}$ and edge set $E(G)$. 
 Denote  by  $\delta(G)$  the \emph{minimum degree} of $G$.
 The \emph{adjacency matrix}
$A(G)$ of $G$  is the $n\times n$ matrix $(a_{ij})$, where
$a_{ij}=1$ if $v_i$ is adjacent to $v_j$, and $0$ otherwise. 
The matrix $Q(G)=D(G)+A(G) $ is known as  the signless Laplacian matrix of $G$, where $D(G)$ is the degree diagonal matrix.
 The \emph{spectral radius} and \emph{signless Laplacian spectral
radius} of $G$ are the largest eigenvalues of $A(G)$ and $Q(G)$, denoted by
$\rho(G)$ and $q(G)$, respectively.

For two disjoint graphs $G$ and $H$,  we denote by $G\bigcup H$ and $G\bigvee H$ the \emph{union} of $G$ and $H$,
 and the \emph{join} of $G$ and $H$ which is  obtained from $G\bigcup H$ by joining every vertex of $G$ to every vertex of $H$, respectively.
Moreover,  $kG$  denotes a graph consisting of $k$ disjoint copies of  $G$.
 Denote by  $G[X,Y]$  the  subgraph of $G$ with all possible edges  with one end vertex in $X$ and the other in $Y$ respectively. Denote
 $e(X,Y)=|E(G[X,Y])|$.

A cycle (path) in a graph $G$ that contains every vertex of $G$ is called a \emph{Hamiltonian cycle (path)} of $G$, respectively.
 A graph $G$ is said to be \emph{Hamiltonian} if it contains a Hamiltonian cycle.
 A bipartite graph $G=(X,Y;E)$ is called \emph {(nearly)} \emph {balanced} if ($|X|-|Y|=1$) $|X|=|Y|$ respectively.
A  (nearly) balanced bipartite graph  $G=(X,Y;E)$ with ($|X|-|Y|=1$) $|X|=|Y|$  is called \emph{Hamilton-biconnected} if for (any two distinct vertices $u,v \in X$) any vertex $u \in X$ and another vertex $v\in Y$, $G$ has a  Hamiltonian path between $u$ and $v$, respectively.
A (nearly) balanced bipartite   graph $G$ is said to be \emph{$2p$-Hamilton-biconnected} if for any balanced subset $W$ of size $2p$ of $V(G)$, the subgraph induced by $V(G)\backslash W$ is  Hamilton-biconnected, respectively. Obviously for $p=0$, $2p$-Hamilton-biconnected graphs are exactly Hamilton-biconnected graphs. For graph notation and terminology undefined here,  readers are referred to \cite{BM}.



 Denote by  $M_{n,m}^{s,t}$ a bipartite graph  obtained from $K_{s,m-t}\bigcup K_{n-s,t}$ by joining
 every vertex in $X_2$ to every vertex in $Y_1$,
 where $K_{s,m-t}=(X_1,Y_1;E_1)$ and $K_{n-s,t}=(X_2,Y_2;E_2)$ with  $|X_1|=s$, $|Y_1|=m-t$, $|X_2|=n-s$, and $|Y_2|=t$ (see Fig.~1).

\medskip
 Denote by $N_{n,n}^{p,1}$  a balanced bipartite graph  obtained from $K_{n-p-2,n-p-2}\bigcup K_{p+1,p+1}\\
 \bigcup K_2$ by joining
 every vertex in $X_1$ to every vertex in $Y_2$, every vertex in $X_2$ to every vertex in $Y_1\bigcup Y_3$, and every vertex in $X_3$ to every vertex in $Y_2$, where $K_{n-p-2,n-p-2}=(X_1,Y_1;E_1)$, $K_{p+1,p+1}=(X_2,Y_2;E_2)$, and  $K_2=(X_3,Y_3;E_3)$ with $|X_1|=|Y_1|=n-p-2$, $|X_2|= |Y_2|=p+1$,
 and $|X_3|= |Y_3|=1$ (see Fig.~1).

\medskip
The problem of deciding whether a  graph is Hamiltonian is NP-complete.  So researchers focus on  giving  reasonable sufficient or necessary conditions  for Hamiltonian cycles in graphs and bipartite graphs.

%
%
%
%
%
%
Moon and Moser \cite{MM} studied balanced bipartite graphs and showed a sufficient   condition for Hamiltonian  cycles in  balanced bipartite graphs with large minimum degree.

\begin{Theorem}\label{Thm1.3}\cite{MM}
Let $G$ be a balanced bipartite graph of order $2n$ with $\delta(G)\geq k$, where $1\leq k\leq\frac{n}{2}$. If
$$e(G)>\max\Bigg\{n(n-k)+k^2,n\bigg(n-\Big\lfloor\frac{n}{2}\Big\rfloor\bigg)+\Big\lfloor\frac{n}{2}\Big\rfloor^2\Bigg\},$$
then $G$ is Hamiltonian.
\end{Theorem}

%
%

 Amar et al. \cite{AFMO} proved a sufficient condition for $2p$-Hamilton-biconnnectedness of balanced bipartite graphs.
 \begin{Theorem}\label{R-Thm1.7}\cite{AFMO}
 Let $p\geq0$ and $G$ be a balanced bipartite graph  of order $2n$. If
 \begin{eqnarray}
   e(G)>n(n-1)+p+1,
 \end{eqnarray}
then $G$ is $2p$-Hamilton-biconnnected.
 \end{Theorem}

%
%
%
%

Recently, Li and Ning \cite{LN} gave the spectral analogue of Moon--Moser's theorem \cite{MM}.
For more results, readers are referred to \cite{Adamus,B,FN,LN1,LLT,N,Z}.

In this paper, we   establish the analogues  of  Moon--Moser's theorem  for $2p$-Hamilton-biconnnectedness of balanced bipartite graphs and nearly balanced bipartite graphs, respectively.


\begin{Theorem}\label{Thm1.7}
Let $p\geq0$ and $G$ be a  balanced bipartite graph of order $2n$ with $\delta(G)\geq k$,  where $n\geq 2k-p+2$.  If $k\geq p$ and
\begin{eqnarray}
  e(G)>n(n-k+p-1)+(k+2)(k-p+1),
\end{eqnarray}
then
$G$ is $2p$-Hamilton-biconnected, unless 
 $G\subseteq M_{n,n}^{n-k,k-p}  $ for  $k\geq p+1$, or $G\subseteq  N_{n,n}^{p,1}$ for $k=p+2$.
\end{Theorem}

\medskip
\noindent {\bf Remark 1.}  Theorem~\ref{R-Thm1.7} \cite{AFMO} and Theorem~\ref{Thm1.7}  are not comparable.
For  $k\geq p+1$ and large $n$,  the condition (2) in Theorem~\ref{Thm1.7} is weaker than  the condition (1) in Theorem~\ref{R-Thm1.7}.

\begin{Theorem}\label{Thm1.8}
Let $p\geq0$ and  $G$ be a  nearly balanced bipartite graph  of order $2n-1$ with   $\delta(G)\geq k$,  where  $n\geq 2k-p+2$. If $k\geq p$ and
$$e(G)> n(n-k+p-2)+(k+2)(k-p+1),$$ then
$G$ is $2p$-Hamilton-biconnected, unless one of the following holds:\\
(i) $G\subseteq M_{n,n-1}^{n-k-1,k-p}$ for $k\geq p+1$;\\
(ii) $G\subseteq M_{n,n-1}^{k-p,n-k-1}$ for $k\geq p+1$;\\
(iii) $G\subseteq M_{n,n-1}^{n-k,k-p-1}$ for $k\geq p+2$.
\end{Theorem}

 Theorems~\ref{Thm1.7} and \ref{Thm1.8}  can be used to obtain some spectral conditions for $2p$-Hamilton-biconnectedness of balanced bipartite  graphs and nearly balanced bipartite graphs in terms of spectral radius or signless Laplacian spectral radius, respectively.

\medskip
For balanced bipartite graphs, we have

\begin{Theorem}\label{Thm1.9}
Let  $p\geq0$, $k\geq p+1$, and $G$ be a  balanced bipartite graph of order $2n$  with $\delta(G)\geq k$. \\
(i) If   $k=p+2$, $n\geq 2k^2+3$, and $\rho(G)\geq \rho(N_{n,n}^{k-2,1})$,  then $G$ is  $2p$-Hamilton-biconnected unless $G= N_{n,n}^{k-2,1}$.\\
(ii) If   $k\neq p+2$,  $n\geq (k+2)(k-p+1)$, and $\rho(G)\geq \rho(M_{n,n}^{n-k,k-p})$,  then $G$ is $2p$-Hamilton-biconnected unless  $G= M_{n,n}^{n-k,k-p}$.
\end{Theorem}

\begin{Theorem}\label{Thm1.10}
Let  $p\geq0$ and $G$ be a  balanced bipartite graph of order $2n$ with  $\delta(G)\geq k $, where  $n\geq (k+2)(k-p+1)$.
If $k\geq p+1$ and  $q(G)\geq q(M_{n,n}^{n-k,k-p})$,  then $G$ is $2p$-Hamilton-biconnected unless  $G= M_{n,n}^{n-k,k-p}$.
\end{Theorem}

For nearly balanced bipartite graphs, we have

\begin{Theorem}\label{Thm1.11}
Let $p\geq0$ and $G$ be a nearly balanced bipartite graph of order $2n-1$ with $\delta(G)\geq k$.  \\
(i) If  $k=p+1$, $n\geq2k+3$, and $\rho(G)\geq \rho(M_{n,n-1}^{1,n-k-1})$,  then $G$ is  $2p$-Hamilton-biconnected unless  $G= M_{n,n-1}^{1,n-k-1}$.\\
(ii) If  $k\geq p+2$, $n\geq \frac{(k+2)(k-p+1)}{2}$,  and $\rho(G)\geq \rho(M_{n,n-1}^{n-k,k-p-1})$,  then $G$ is $2p$-Hamilton-biconnected unless  $G= M_{n,n-1}^{n-k,k-p-1}$.
\end{Theorem}

\begin{Theorem}\label{Thm1.12}
Let $p\geq0$ and $G$ be a nearly balanced bipartite graph of order $2n-1$ with $\delta(G)\geq k $.  \\
(i) If  $k=p+1$, $n\geq2k+4$, and $q(G)\geq q(M_{n,n-1}^{n-k-1,1})$,  then $G$ is  $2p$-Hamilton-biconnected unless  $G= M_{n,n-1}^{n-k-1,1}$.\\
(ii) If  $k\geq p+2$, $n\geq \frac{(k+2)(k-p+1)}{2}$, and $q(G)\geq q(M_{n,n-1}^{n-k,k-p-1})$,  then $G$ is $2p$-Hamilton-biconnected unless  $G= M_{n,n-1}^{n-k,k-p-1}$.
\end{Theorem}

The rest of this paper is organized as follows. In Section~2, we state some known and  new results that will be
used in the proofs of Theorems~\ref{Thm1.7}--\ref{Thm1.12}.
In Section~3, we present some necessary lemmas and prove Theorems~\ref{Thm1.7} and  \ref{Thm1.8}.
In Section~4, we present some necessary lemmas and prove Theorems~\ref{Thm1.9} and \ref{Thm1.10}. Some corollaries are also included.
In Section~5, we present some necessary lemmas and prove Theorems~\ref{Thm1.11} and \ref{Thm1.12}. Some corollaries are also included.

\section{Preliminarily }
Next we introduce some more  terminologies and notations, which will be used in this section and the proofs of Theorems~\ref{Thm1.7} and \ref{Thm1.8}.

Recall that the \emph{k-biclosure} of a bipartite graph $G=(X,Y;E)$ \cite{BC} is the unique smallest bipartite graph $H$ of order $|V(H)|:=|V(G)|$ such that $G\subseteq H$ and $d_H(x)+d_H(y)<k$ for any two  non-adjacent vertices $x\in X$ and $y\in Y$. The $k$-biclosure of $G$ is denoted by $cl_k(G)$, and can be obtained from $G$ by a recursive procedure which consists of joining non-adjacent vertices in different classes with degree sum at least $k$ until no such pair remains. A bipartite graph is called \emph{$k$-closed} if $G=cl_k(G)$.


\bigskip

\begin{tikzpicture}[scale=1.6]
\path (-1,0)  coordinate (P1);\path (-0.3,0)  coordinate (P2);  \path (0.3,0)  coordinate (P3); \path (1,0)  coordinate (P4);
\path (-1,-1)  coordinate (Q1);\path (-0.3,-1)  coordinate (Q2);  \path (0.3,-1)  coordinate (Q3); \path (1,-1)  coordinate (Q4);
\path (-0.8,0)  coordinate (R1);\path (-0.65,0)  coordinate (R2);  \path (-0.5,0)  coordinate (R3);
\path (0.5,0)  coordinate (R4); \path (0.65,0)  coordinate (R5);\path (0.8,0)  coordinate (R6);
\path (-0.8,-1)  coordinate (R7);\path (-0.65,-1)  coordinate (R8);  \path (-0.5,-1)  coordinate (R9);
\path (0.5,-1)  coordinate (R10); \path (0.65,-1)  coordinate (R11);\path (0.8,-1)  coordinate (R12);
 \foreach \i in {1,2,3,4}
{\fill (P\i) circle (1.5pt);\fill (Q\i) circle (1.5pt);}
\foreach \i in {1,...,12}
{\fill (R\i) circle (0.5pt);}
 \foreach \i in {1,2}
 { \foreach \j in {1,2,3,4}
{\draw (Q\i)--(P\j);}}
 \foreach \i in {3,4}
 { \foreach \j in {3,4}
{\draw (P\i)--(Q\j);}}
\path (-0.65,0.4) node (s) {$s$};
\path (-0.65,0.2) node (){$\overbrace{~~~~~~~}$};
\path (0.65,0.4) node (n-s) {$n-s$};
\path (0.65,0.2) node (){$\overbrace{~~~~~~~}$};
\path (-0.65,-1.4) node (m-t) {$m-t$};
\path (-0.65,-1.2) node () {$\underbrace{~~~~~~~}$};
\path (0.65,-1.4) node () {$t$};
\path (0.65,-1.2) node () {$\underbrace{~~~~~~~}$};
\path (0,-1.8) node () {$M_{n,m}^{s,t}$};

\path (1.5,0)  coordinate (P1);\path (2.2,0)  coordinate (P2);  \path (2.8,0)  coordinate (P3); \path (3.5,0)  coordinate (P4);
\path (1.5,-1)  coordinate (Q1);\path (2.2,-1)  coordinate (Q2);  \path (2.8,-1)  coordinate (Q3); \path (3.5,-1)  coordinate (Q4);
\path (2.5,0)  coordinate (S1);\path (2.5,-1)  coordinate (S2);
\path (1.7,0)  coordinate (R1);\path (1.85,0)  coordinate (R2);  \path (2,0)  coordinate (R3);
\path (3,0)  coordinate (R4); \path (3.15,0)  coordinate (R5);\path (3.3,0)  coordinate (R6);
\path (1.7,-1)  coordinate (R7);\path (1.85,-1)  coordinate (R8);  \path (2,-1)  coordinate (R9);
\path (3,-1)  coordinate (R10); \path (3.15,-1)  coordinate (R11);\path (3.3,-1)  coordinate (R12);
\fill (S1) circle (1.5pt);\fill (S2) circle (1.5pt);
 \foreach \i in {1,2,3,4}
{\fill (P\i) circle (1.5pt);\fill (Q\i) circle (1.5pt);}
\foreach \i in {1,...,12}
{\fill (R\i) circle (0.5pt);}
 \foreach \i in {1,2}
 { \draw (S1)--(Q\i);\foreach \j in {1,2,3,4}
{\draw (Q\i)--(P\j);}}
 \foreach \i in {3,4}
 {\draw (S2)--(P\i); \foreach \j in {3,4}
{\draw (P\i)--(Q\j);}}
\draw (S1)--(S2);
\path (1.85,0.4) node () {$s-1$};
\path (1.85,0.2) node (){$\overbrace{~~~~~~~}$};
\path (3.15,0.4) node () {$n-s$};
\path (3.15,0.2) node (){$\overbrace{~~~~~~~}$};
\path (1.8,-1.4) node () {$m-t-1$};
\path (1.8,-1.2) node () {$\underbrace{~~~~~~~}$};
\path (3.15,-1.4) node () {$t$};
\path (3.15,-1.2) node () {$\underbrace{~~~~~~~}$};
\path (2.5,-1.8) node () {$M_{n,m}^{s,t;-}$};

\path (4,0)  coordinate (P1);\path (4.5,0)  coordinate (P2); \path (4.9,0)  coordinate (P3);\path (5.4,0)  coordinate (P4);
 \path (6,0)  coordinate (P5); \path (6.5,0)  coordinate (P6); \path (6.7,0)  coordinate (P7); \path (7.2,0)  coordinate (P8);
\path (4.4,-1)  coordinate (Q1);\path (5.2,-1)  coordinate (Q2);  \path (5.9,-1)  coordinate (Q3); \path (6.7,-1)  coordinate (Q4);
\path (4.15,0)  coordinate (R1);\path (4.25,0)  coordinate (R2);  \path (4.35,0)  coordinate (R3);
\path (5.05,0)  coordinate (R4); \path (5.15,0)  coordinate (R5);\path (5.25,0)  coordinate (R6);
\path (6.15,0)  coordinate (R7);\path (6.25,0)  coordinate (R8);  \path (6.35,0)  coordinate (R9);
\path (6.85,0)  coordinate (R10); \path (6.95,0)  coordinate (R11);\path (7.05,0)  coordinate (R12);
\path (4.6,-1)  coordinate (R13); \path (4.8,-1)  coordinate (R14);\path (5,-1)  coordinate (R15);
\path (6.1,-1)  coordinate (R16); \path (6.3,-1)  coordinate (R17);\path (6.5,-1)  coordinate (R18);
\fill (S1) circle (1.5pt);\fill (S2) circle (1.5pt);
 \foreach \i in {1,...,8}
{\fill (P\i) circle (1.5pt);}
 \foreach \i in {1,...,4}
{\fill (Q\i) circle (1.5pt);}
\foreach \i in {1,...,18}
{\fill (R\i) circle (0.5pt);}
 \foreach \i in {1,2}
 {\foreach \j in {1,...,8}
{\draw (Q\i)--(P\j);}}
 \foreach \i in {3,4}
 { \foreach \j in {3,4}
{\draw (P\i)--(Q\j);}}
\draw (Q3)--(P5);\draw (Q3)--(P6);\draw (Q4)--(P7);\draw (Q4)--(P8);
\path (4.3,0.4) node () {$t$};
\path (4.3,0.2) node () {$\overbrace{~~~~~}$};
\path (5.2,0.4) node () {$l$};
\path (5.2,0.2) node () {$\overbrace{~~~}$};
\path (6.2,0.4) node () {$k-l$};
\path (6.2,0.2) node () {$\overbrace{~~~}$};
\path (7,0.4) node () {$k-l$};
\path (7,0.2) node (){$\overbrace{~~~}$};
\path (4.8,-1.4) node () {$m-k+p$};
\path (4.8,-1.2) node (){$\underbrace{~~~~~~~}$};
\path (6.3,-1.4) node () {$k-p$};
\path (6.3,-1.2) node (){$\underbrace{~~~~~~~}$};
\path (5.5,-1.8) node () {$F_{n,m}^{k,p,l}$, };
\path (5.5,-2.1) node ()  { $t=n-(k-p)(k-l)-l$};
\end{tikzpicture}

\begin{tikzpicture}[scale=1.55]
\path (1.5,0)  coordinate (P1);\path (2.2,0)  coordinate (P2);  \path (2.6,0)  coordinate (P3);\path (3.3,0)  coordinate (P4);
 \path (3.7,0)  coordinate (P5);\path (3.7,-1)  coordinate (P6);
\path (1.5,-1)  coordinate (Q1);\path (2.2,-1)  coordinate (Q2);  \path (2.6,-1)  coordinate (Q3); \path (3.3,-1)  coordinate (Q4);
\path (1.7,0)  coordinate (R1);\path (1.85,0)  coordinate (R2);  \path (2,0)  coordinate (R3);
\path (2.8,0)  coordinate (R4); \path (2.95,0)  coordinate (R5);\path (3.1,0)  coordinate (R6);
\path (1.7,-1)  coordinate (R7);\path (1.85,-1)  coordinate (R8);  \path (2,-1)  coordinate (R9);
\path (2.8,-1)  coordinate (R10); \path (2.95,-1)  coordinate (R11);\path (3.1,-1)  coordinate (R12);
\fill (P5) circle (1.5pt);\fill (P6) circle (1.5pt);
 \foreach \i in {1,2,3,4}
{\fill (P\i) circle (1.5pt);\fill (Q\i) circle (1.5pt);}
\foreach \i in {1,...,12}
{\fill (R\i) circle (0.5pt);}
 \foreach \i in {1,2,3,4}
 { \foreach \j in {1,2,3,4}
{\draw (Q\i)--(P\j);}}
 \foreach \i in {3,4}
 { \foreach \j in {3,4}
{\draw (P\i)--(Q\j);}}
\draw (P5)--(Q3);\draw (P5)--(Q4);\draw (P6)--(P3);\draw (P6)--(P4);\draw (P6)--(P4);\draw (P5)--(P6);
\path (1.8,0.4) node () {$n-p-2$};
\path (1.8,0.2) node () {$\overbrace{~~~~~~~}$};
\path (3,0.4) node () {$p+1$};
\path (3,0.2) node () {$\overbrace{~~~~~~~}$};
\path (1.8,-1.4) node () {$n-p-2$};
\path (1.8,-1.2) node () {$\underbrace{~~~~~~~}$};
\path (3,-1.4) node () {$p+1$};
\path (3,-1.2) node () {$\underbrace{~~~~~~}$};
\path (2.5,-1.8) node () {$N_{n,n}^{p,1}$};

\path (5.5,0)  coordinate (P1);\path (6.2,0)  coordinate (P2);  \path (6.6,0)  coordinate (P3);\path (7.3,0)  coordinate (P4);
 \path (7.7,0)  coordinate (P5);\path (7.7,-1)  coordinate (P6);
\path (5.5,-1)  coordinate (Q1);\path (6.2,-1)  coordinate (Q2);  \path (6.6,-1)  coordinate (Q3); \path (7.3,-1)  coordinate (Q4);
\path (5.7,0)  coordinate (R1);\path (5.85,0)  coordinate (R2);  \path (6,0)  coordinate (R3);
\path (6.8,0)  coordinate (R4); \path (6.95,0)  coordinate (R5);\path (7.1,0)  coordinate (R6);
\path (5.7,-1)  coordinate (R7);\path (5.85,-1)  coordinate (R8);  \path (6,-1)  coordinate (R9);
\path (6.8,-1)  coordinate (R10); \path (6.95,-1)  coordinate (R11);\path (7.1,-1)  coordinate (R12);
\fill (P5) circle (1.5pt);\fill (P6) circle (1.5pt);
 \foreach \i in {1,2,3,4}
{\fill (P\i) circle (1.5pt);\fill (Q\i) circle (1.5pt);}
\foreach \i in {1,...,12}
{\fill (R\i) circle (0.5pt);}
 \foreach \i in {1,2,3,4}
 { \foreach \j in {1,2,3,4}
{\draw (Q\i)--(P\j);}}
 \foreach \i in {3,4}
 { \foreach \j in {3,4}
{\draw (P\i)--(Q\j);}}
\draw (P5)--(Q3);\draw (P5)--(Q4);\draw (P6)--(P3);\draw (P6)--(P4);
\path (5.8,0.4) node () {$n-p-3$};
\path (5.8,0.2) node () {$\overbrace{~~~~~~}$};
\path (7,0.4) node () {$p+2$};
\path (7,0.2) node () {$\overbrace{~~~~~~~}$};
\path (5.8,-1.4) node () {$n-p-3$};
\path (5.8,-1.2) node () {$\underbrace{~~~~~~~}$};
\path (7,-1.4) node () {$p+2$};
\path (7,-1.2) node () {$\underbrace{~~~~~~~~}$};
\path (6.5,-1.8) node () {$N_{n,n}^{p,2}$};
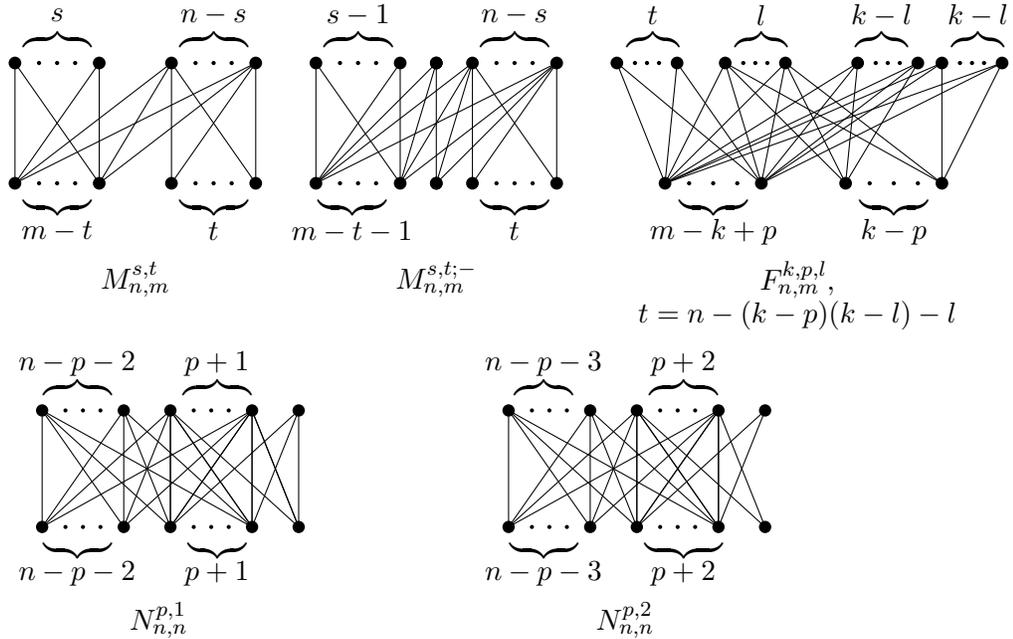
\path (5,-2.4) node () {Fig.~1. Graphs $M_{n,m}^{s,t}$, $M_{n,m}^{s,t;-}$, $F_{n,m}^{k,p,l}$, $N_{n,n}^{p,1}$ and $N_{n,n}^{p,2}$.};
\end{tikzpicture}

 Denote by  $M_{n,m}^{s,t;-}$ a bipartite graph  obtained from $K_{s-1,m-t-1}\bigcup K_2\bigcup K_{n-s,t}$ by joining
every vertex in $X_2$ to every vertex in $Y_1$, and every vertex in $Y_3$ to every vertex in $Y_1\bigcup Y_2$,
 where $K_{s-1,m-t-1}=(X_1,Y_1;E_1)$, $K_2=(X_2;Y_2;E_2)$,  and $K_{n-s,t}=(X_3,Y_3;E_3)$ with  $|X_1|=s-1$, $|Y_1|=m-t-1$, $|X_2|=|Y_2|=1$, $|X_3|=n-s$, and $|Y_3|=t$ (see Fig.~1). Obviously $M_{n,m}^{s,t;-}\subseteq M_{n,m}^{s,t}$.

\medskip
Denote by $N_{n,n}^{p,2}$  a balanced bipartite graph  obtained from $K_{n-p-3,n-p-3}\bigcup K_{p+2,p+2}\\
\bigcup \overline{K}_2$ by joining
 every vertex in $X_1$ to every vertex in $Y_2$, every vertex in $X_2$ to every vertex in $Y_1\bigcup Y_3$, and every vertex in $X_3$ to every vertex in $Y_2$, where $K_{n-p-3,n-p-3}=(X_1,Y_1;E_1)$, $K_{p+2,p+2}=(X_2,Y_2;E_2)$, and  $\overline{K}_2=(X_3,Y_3;E_3)$ with $|X_1|=|Y_1|=n-p-3$, $|X_2|=|Y_2|=p+2$,
 and $|X_3|=|Y_3|=1$  (see Fig.~1).

 \medskip

Given integers $n,m,k,p,l$, where  $k\geq p+2$, $0\leq l\leq k-1$, $n\geq(k-p)(k-l)+l$, and  $n-1\leq m\leq n$,
we denote by $F_{n,m}^{k,p,l}$  a  bipartite graph  obtained from $M_{n-(k-p)(k-l)+l,m}^{n-(k-p)(k-l),k-p}$
by attaching $k-l$ pendant vertices at every vertex of those $k-p$ vertices with degree $l$, respectively, and then joining every pendant vertex to
every vertex with degree $n-(k-p)(k-l)$ in $M_{n-(k-p)(k-l)+l,m}^{n-(k-p)(k-l),k-p}$ (see Fig.~1).

\bigskip

The following lemma follows from the Perron--Frobenius theorem.
\begin{Lemma}\label{Lem2.1}
Let $G$ be a connected graph. If $H$ is a (proper) subgraph  of $G$, then $ \rho(H) (<)\leq \rho(G)$
and $ q(H) (<)\leq q(G)$, respectively.
\end{Lemma}

\begin{Lemma}\label{Lem2.2}\cite{BFP}
Let $G$ be a bipartite graph. Then \\
$$\rho(G)\leq \sqrt{e(G)},$$
with equality if and only if $G$ is a disjoint union of a complete bipartite graph and isolated vertices.
\end{Lemma}

\begin{Lemma}\label{Lem2.3}\cite{LN}
Let $G$ be a balanced bipartite graph of order $2n$. Then\\
$$q(G) \leq \frac{e(G)}{n}+n.$$
with equality if and only if  $G= K_{n,n}$.
\end{Lemma}

\medskip

\noindent{\bf Remark~2:} The extremal graph in Lemma~\ref{Lem2.3} is not characterized in \cite{LN}.  But it is easy to obtain the extremal graph by combining the proof of Lemma~\ref{Lem2.3} and Das's bound \cite[Theorem~4.5]{Das}.

\medskip
Note that $G\subseteq  cl_{n+p+1} (G)$. If $G$ is $2p$-Hamilton-biconnected then so is $cl_{n+p+1} (G)$.  Combining this with \cite[Theorem~3.3.1]{AFMO},
we have the following lemma.

\begin{Lemma}\label{Lem3.1}
Let $p\geq 0$ and $G$ be a balanced bipartite graph of order $2n$. Then  $G$ is $2p$-Hamilton-biconnected if and only if $cl_{n+p+2} (G)$ is $2p$-Hamilton-biconnected.
\end{Lemma}


\begin{Lemma}\label{Lem4.2}
Let  $p\geq 0$ and $G$ be a nearly balanced bipartite graph of order $2n-1$. Then  $G$ is $2p$-Hamilton-biconnected if and only if $cl_{n+p+1} (G)$ is $2p$-Hamilton-biconnected.
\end{Lemma}

\begin{Proof}
Since $G\subseteq  cl_{n+p+1} (G)$, if $G$ is $2p$-Hamilton-biconnected then so is $cl_{n+p+1} (G)$. Conversely, suppose that $cl_{n+p+1} (G)$ is  $2p$-Hamilton-biconnected. Denote $G=(X,Y;E) $  with $|X|=n$ and  $|Y|=n-1$.
We show  that if $G+xy$ is $2p$-Hamilton-biconnected for two non-adjacent vertices $x\in X$ and  $y\in Y$ with $d_G(x)+d_G(y)\geq n+p+1$, then $G$
 is $2p$-Hamilton-biconnected. Indeed, if  $G$ is not $2p$-Hamilton-biconnected, then there  exists a balanced subset $W$ of size $2p$ of $V(G)$ and  two vertices $x_1,x_2\in X\backslash W$ such that the subgraph $F$ induced by $V(G)\backslash W$ has no Hamiltonian path between $x_1$ and $x_2$.
On the other hand, since $G+xy$  is $2p$-Hamilton-biconnected, the graph $F+xy$ has a Hamiltonian path between $x_1$ and $x_2$  and thus $x\in X\backslash W$ and $y\in Y\backslash W$. Let $H$ be a graph obtained from $F$ by adding a new vertex $v$ in $Y$ and two edges $vx_1$ and $vx_2$. Then $H$ is not Hamiltonian, but $H+xy$ is Hamiltonian. Note that
\begin{eqnarray*}
  d_H(x)+d_H(y)&\geq& d_F(x)+d_F(y) \geq(d_G(x)-p)+(d_G(y)-p)\\
   &=& d_G(x)+d_G(y)-2p\geq n-p+1=\frac{1}{2}|V(H)|+1.
\end{eqnarray*}
It follows from \cite[Theorem~6.2]{BC} that $H$ is Hamiltonian, a contradiction.
Note that $cl_{n+p+1} (G)$ is a graph obtained from $G$ by a recursive procedure joining non-adjacent vertices in different classes with degree sum at least $n+p+1$ until no such pair remains. Since $cl_{n+p+1} (G)$ is $2p$-Hamilton-biconnected,  $G$ is also $2p$-Hamilton-biconnected.
\end{Proof}

\medskip

The proofs of Lemmas~\ref{R-Lem2.3}--\ref{R3-Lem2.3}are put in the appendix, since they are   technical and complicated.

\begin{Lemma}\label{R-Lem2.3}
$F_{n,m}^{k,p,l}$ is $2p$-Hamilton-biconnected.
\end{Lemma}

\begin{Lemma}\label{R2-Lem2.3}
(i) For $p\geq0$, $s\geq2$, $t\geq1$, and $n=s+t+p+1$,  $M_{n,n}^{s,t;-}$ is $2p$-Hamilton-biconnected.\\
(ii) For $p\geq0$ and $n\geq p+6$, $N_{n,n}^{p,2}$ is $2p$-Hamilton-biconnected.
\end{Lemma}

\begin{Lemma}\label{R3-Lem2.3}
(i) For   $p\geq 0$, $s,t\geq1$, and $ \max \{s+p+2,t+p+2\} \leq n\leq s+t+p+1$, $M_{n,n-1}^{s,t}$ is not  $2p$-Hamilton-biconnected.\\
(ii) For  $p\geq 0$, $s,t\geq1$,  and  $n=s+t+p$, $M_{n,n}^{s,t}$ is not  $2p$-Hamilton-biconnected.\\
(iii) For $p\geq 0$ and $n\geq p+6$, $N_{n,n}^{p,1}$ is not  $2p$-Hamilton-biconnected.
\end{Lemma}

\section{Proofs of Theorems~\ref{Thm1.7} and \ref{Thm1.8}}
In order to prove Theorems~\ref{Thm1.7} and \ref{Thm1.8}, we first prove the following lemma, in which the techniques  are from \cite[Lemma~4]{LN}.

\begin{Lemma}\label{Lem2.4}
Let $G$ be an $(n+p+2)$-closed balanced bipartite graph  of order $2n$, where  $k\geq p\geq0$ and $n\geq 2k-p+2$. If
$$e(G)> n(n-k+p-1)+(k+2)(k-p+1),$$
then $G$ contains a complete
bipartite graph of order $2n-k+p$. Furthermore, if $\delta(G)\geq k$,
then $K_{n,n-k+p}\subseteq G$, or $G\in \{N_{n,n}^{p,1},N_{n,n}^{p,2}\}$ for $k=p+2$.
\end{Lemma}

\begin{Proof}
Denote $G=(X,Y;E)$ with $|X|=|Y|=n$. Let  $U=\big\{x\in X: d_G(x)\geq\frac{n+p+2}{2}\big\}$ and  $W=\big\{y\in Y: d_G(y)\geq\frac{n+p+2}{2}\big\}$. Then
$$n(n-k+p-1)+(k+2)(k-p+1)< e(G)\leq n|U|+\frac{(n-|U|)(n+p+1)}{2}.$$
Since $k\geq p$ and $n\geq 2k-p+2$, we have
$$|U|\geq\frac{n^2 - (2k-p+3)n + 2(k + 2)(k - p + 1) + 2}{n-p-1}> k+1,$$ which implies that $|U|\geq k+2$.
By symmetry, $|W|\geq k+2$.
Since $G$ is an $(n+p+2)$-closed balanced bipartite graph,
every vertex in $U$ is adjacent to every vertex in $W$ and thus $K_{k+2,k+2}\subseteq G$.
Let $t$ be the largest integer such that $K_{t,t}\subseteq G$.

\medskip

{\bf{Claim~1}.} $t\geq n-k+p$.

\medskip

Suppose that $k+2\leq t\leq n-k+p-1$. Let $X_1\subseteq X$ and $Y_1\subseteq Y$ with $|X_1|=|Y_1|=t$ such that $G[X_1, Y_1]=K_{t,t}$.
Set $X_2=X\backslash X_1$ and $Y_2=Y\backslash Y_1$. Since $t$ is  the largest integer such that $K_{t,t}\subseteq G$,
 there exists a corresponding vertex $y\in Y_1$ such that $xy\notin E(G)$ for every  $x\in X_2$  (by symmetry). It follows that  $d_G(x)\leq n+p-t+1$ for every  $x\in X_2$.
Hence
\begin{eqnarray*}
  e(G) &=& e(X_1, Y_1)+e(X_1, Y_2)+e(X_2,Y) \\
   &\leq& t^2+t(n-t)+(n+p-t+1)(n-t)\\
   &=& t^2-(n+p+1)t+n(n+p+1)\\
   &\leq&(n-k+p-1)^2-(n+p+1)(n-k+p-1)+n(n+p+1)\\
   &=& n(n-k+p-1)+(k+2)(k-p+1)\\
   &<&e(G),
\end{eqnarray*}
a contradiction. Thus Claim~1 holds.

Let $s$ be the largest integer  such that $K_{s,t}\subseteq G$.  Obviously, $s\geq t$.
Let $X_1\subseteq X$ and $Y_1\subseteq Y$ such that $G[X_1, Y_1]=K_{s,t}$, where $|X_1|=s$ and  $ |Y_1|=t$.
Set $X_2=X\backslash X_1$ and $Y_2=Y\backslash Y_1$.

\medskip

{\bf{Claim~2}.} $s+t\geq 2n-k+p$.

\medskip

Suppose that $s+t\leq  2n-k+p-1$.
It follows from Claim~1 that $n-k+p\leq t\leq n-\frac{k-p+1}{2}$ and $t\leq s\leq 2n-k+p-t-1$.
Since  $G$ is an $(n+p+2)$-closed balanced bipartite graph,  $d_G(x)\leq n+p-s+1$ for every $x\in X_2$ and  $d_G(y)\leq n+p-t+1$ for every $y\in Y_2$.
Hence
\begin{eqnarray*}
  e(G) &\leq& e(X_1, Y_1)+e(X_2, Y)+e(X, Y_2)\\
   &\leq&st+(n+p-s+1)(n-s)+(n+p-t+1)(n-t)\\
   &=& s^2 - (2n+ p - t  +1)s  + (n - t)(n + p + t +1)+ n(n + p + 1)\\
   &\leq& (2n-k+p-t-1)^2-(2n+ p  - t  +1)(2n-k+p-t-1)+\\
   && (n - t)(n + p  - t +1)+n(n + p + 1)\\
   &=&t^2 - (2n -k+ p  -  1)t +2n^2 - 2(k -p+1)n+(k + 2)(k - p +1)\\
   &\leq& (n-k+p)^2-(2n -k+ p -1)(n-k+p)+2n^2 +2(k -p+1)n+\\
   &&(k + 2)(k - p +1)\\
   &=&n(n-k+p-1)+(k + 1)(k - p + 1)+1\\
   &< & e(G),
\end{eqnarray*}
a contradiction. Thus Claim~2 holds.

\medskip
It follows  from Claim~2 that  $K_{s,t}$ is a complete bipartite graph of order at least $2n-k+p$. Hence $G$ contains a complete bipartite graph of order $2n-k+p$.

\medskip

{\bf{Claim~3}.} If $\delta(G)\geq k$, then
$K_{n,n-k+p}\subseteq G$, or $G\in \{N_{n,n}^{p,1},N_{n,n}^{p,2}\}$ for $k=p+2$.

\medskip

If $t=n-k+p$, then  Claim~2 implies that $s=n$ and thus $K_{n,n-k+p}\subseteq G$.
So we can assume that $t\geq n-k+p+1$. Next we consider the following two cases.

\medskip

{\bf{Case~1}.} $s=n-k+p+1$. Obviously, $k\geq p+1$ and $t=n-k+p+1$.
If $k\geq p+3$, then $s+t=2(n-k+p+1)<2n-k+p$, which contradicts  Claim~2.
Hence $ p+1\leq k\leq p+2$. Furthermore, if $k=p+1$, then $G=K_{n,n}$ and thus $K_{n,n-k+p}\subseteq G$.
If $k=p+2$, then $s=t=n-1$.
 Note that $G$ is an $(n+p+2)$-closed bipartite graph with $\delta(G)\geq k$.
  If there exists a vertex  $v\in V(G)\backslash V(K_{n-1,n-1})$ such that $d_G(v)\geq k+1$, then $K_{n,n-1}\subseteq G$.
If $d_G(v)=k=p+2$  for every $v\in V(G)\backslash V(K_{n-1,n-1})$,  then $G= N_{n,n}^{p,1}$  or   $G= N_{n,n}^{p,2}$ for $k=p+2$.
\medskip

{\bf{Case~2}.} $s\geq n-k+p+2$.  Clearly, $d_G(y)\geq n-k+p+2$ for every $y\in Y_1$. Then  every vertex in $Y_1$ is adjacent to every vertex in $X$.  This implies that $K_{n,n-k+p}\subseteq G$.
\end{Proof}

\begin{Corollary}\label{Cor2.5}
Let $G$ be an $(n+p+1)$-closed nearly balanced bipartite graph of order $2n-1$, where   $k\geq p\geq 0$ and $n\geq 2k-p+2$. If
 $$e(G)>n(n-k+p-2)+(k+2)(k-p+1),$$ then $G$ contains a complete
bipartite graph of order $2n-k+p-1$. Furthermore, if $\delta(G)\geq k$, then
 $K_{n-1,n-k+p}\subseteq G$ or $K_{n,n-k+p-1}\subseteq G$.
\end{Corollary}

\begin{Proof}
Denote $G=(X,Y;E)$ with $|X|=n$ and $|Y|=n-1$. Let
 $H$ be a graph with vertex set $V(G)\bigcup \{y\}$ and edge set $E(G)\bigcup \{xy:x\in X\}$, where $y\notin V(G)$. Then $H$ is an  $(n+p+2)$-closed balanced bipartite graph of order $2n$ and $e(H)= e(G)+n> n(n-k+p-1)+(k+2)(k-p-1)$. By Lemma~\ref{Lem2.4}, $H$ contains a complete bipartite graph of order $2n-k+p$. Thus $G$ contains a complete bipartite graph of order $2n-k+p-1$.
Note that if $\delta(G)\geq k$, then  $\delta(H)\geq \delta(G)\geq k$.  It follows from Lemma~\ref{Lem2.4} that $K_{n,n-k+p}\subseteq H$, or $H\in \{N_{n,n}^{p,1},N_{n,n}^{p,2}\}$ for $k=p+2$. If  $K_{n,n-k+p}\subseteq H$,
then $K_{n,n-k+p-1}\subseteq G$ or $K_{n-1,n-k+p}\subseteq G$. If $H\in\{N_{n,n}^{p,1},N_{n,n}^{p,2}\}$ for $k=p+2$, then $\delta(G)=p+1<k$, a contradiction.
This completes the proof.
\end{Proof}

\begin{Lemma}\label{Lem2.6}
Let $p\geq0$ and  $G=(X,Y;E)$ be an $(m+p+2)$-closed  bipartite graph with $|X|=n$, $|Y|=m$, and  $\delta(G)\geq k$, where
  $n \geq 2k-p+2$ and $n-1\leq m\leq n$. Suppose that $k\geq  p$ and $K_{n,m-k+p}\subseteq G$.\\
 (i)   If $m=n$, then $G$ is $2p$-Hamilton-biconnected, unless $G= M_{n,n}^{n-k,k-p}$ for $k\geq p+1$.\\
(ii)  If  $m=n-1$, then  $G$ is $2p$-Hamilton-biconnected, unless one of the following holds:

(a) $G= M_{n,n-1}^{n-k,k-p}$ for $k\geq p+1$;

(b) $G= M_{n,n-1}^{n-k-1,k-p}$ for $k\geq p+1$;

(c) $G= M_{n,n-1}^{n-k,k-p-1;-} $ for $k\geq p+2$;

(d) $G= M_{n,n-1}^{n-k,k-p-1}$ for $k\geq p+2$.
\end{Lemma}

\begin{Proof}
Suppose that $G$ is not $2p$-Hamilton-biconnected.
Let $t$ be the largest integer such that $K_{n,t}\subseteq G$, and
$Y_1\subseteq Y$ such that $G[X,Y_1]=K_{n,t}$. Obviously, $t\geq m-k+p$.
We claim that $m-k+p\leq t\leq m-k+p+1$. Note that  $G$ is an $(m+p+2)$-closed  bipartite graph and $\delta(G)\geq k$. If $t> m-k+p+1$, then every vertex in $Y$ is adjacent to every vertex in $X$, and thus $G=K_{n,m}$, a contradiction.
Next we consider the following two cases.

\medskip
{\bf{Case~1}.} $t=m-k+p$.  Then $|Y_1|=m-k+p$ and $|Y\backslash Y_1|=k-p\geq1$.
We show that $k\leq d_G(y)\leq k+1$  for every $y\in Y\backslash Y_1$.   Indeed, if there exists a vertex  $y\in Y\backslash Y_1$ such that $d_G(y)\geq k+2$,
then $y$ is adjacent to every vertex in $X$ and thus $t\geq m-k+p+1$, a contradiction.  Next we consider the following two subcases.

\medskip
{\bf{Case~1.1}.} For every $y\in Y\backslash Y_1$, $d_G(y)=k$.
Set $X=\bigcup _{i=1}^3 X_i$, where  $X_1=\{x\in X:d_G(x)=m-k+p\}$,  $X_2=\{x\in X: d_G(x)=m-k+p+1\}$, and $X_3=\{x\in X: d_G(x)\geq m-k+p+2\}$.
Set $Y_2=Y\backslash Y_1$ and $l=|X_3|$.
Since $G$ is an $(m+p+2)$-closed  bipartite graph with $\delta(G)\geq k$, every vertex in $Y_2$ is adjacent to every vertex in $X_3$. This implies that $0\leq l\leq k$.
Furthermore,  every  vertex in $Y_2$ is adjacent to $k-l$ vertices in $X_2$ and any two distinct vertices in $Y_2$ have no common neighbors in $X_2$.
This implies that $|X_2|=(k-p)(k-l)$.
Moreover, if $k\geq p+2$ and  $0\leq l\leq k-1$, then $G= F_{n,m}^{k,p,l} $. By Lemma~\ref{R-Lem2.3}, $F_{n,m}^{k,p,l}$ is $2p$-Hamilton-biconnected, a contradiction. If  $l=k\geq p+2$ or $k=p+1$, then $G= M_{n,m}^{n-k,k-p}$. It follows from Lemma~\ref{R3-Lem2.3}~($i$) and ($ii$) that $M_{n,m}^{n-k,k-p}$ is not $2p$-Hamilton-biconnected, as desired.

\medskip
{\bf{Case~1.2}.}  There exists a vertex $y\in Y\backslash Y_1$ such that $ d_G(y)=k+1$.
Set $X=\bigcup _{i=1}^2 X_i$ and $Y=\bigcup _{i=1}^3 Y_i$, where $X_1=\{x\in X:d_G(x)=m-k+p\}$, $X_2=\{x\in X: d_G(x)\geq m-k+p+1\}$, $Y_2=\{y\in Y:d_G(y)=k\}$, and $Y_3=\{y\in Y: d_G(y)=k+1\}$. Since  $G$ is an $(m+p+2)$-closed  bipartite graph with $\delta(G)\geq k$,  every vertex in $Y_3$ is adjacent to every vertex in $X_2$. This implies that $|X_2|=k+1$ and thus $|X_1|=n-k-1$.

We first assume that $Y_2= \emptyset$. It is easy to see that $G= M_{n,m}^{n-k-1,k-p}$ for $k\geq p+1$.
Suppose that $m=n$.  Since  $ M_{n,n}^{n-k-1,k-p;-}$ is a spanning subgraph of $ M_{n,n}^{n-k-1,k-p}$, it follows from   Lemma~\ref{R2-Lem2.3}~($i$) that
$ M_{n,n}^{n-k-1,k-p}$ is $2p$-Hamilton-biconnected, a contradiction. Next suppose that  $m=n-1$. By Lemma~\ref{R3-Lem2.3}~($i$), $M_{n,n-1}^{n-k-1,k-p}$ is not $2p$-Hamilton-biconnected, as desired.

We next assume that $Y_2\neq\emptyset$.  We show that $|Y_3|=1$. Indeed, if $|Y_3|\geq2$, then $d_G(x)\geq m-k+p+2$ for every $x\in X_2$ and hence every vertex in $Y_2$ is adjacent to every vertex in $X_2$. This implies that $d_G(y)= k+1$ for every $y\in Y_2$, a contradiction.
By a similar argument to the proof of $|Y_3|=1$,  there exists a vertex $x\in X_2$ adjacent to none of vertices in $Y_2$.
Moreover, since $\delta(G)\geq k$, every vertex in $Y_2$ is adjacent to every vertex in $X_2\backslash \{x\}$.
Hence $G= M_{n,m}^{n-k,k-p-1;-} $ for $k\geq p+2$.
Suppose that $m=n$.  By Lemma~\ref{R2-Lem2.3}~($i$), $M_{n,n}^{n-k,k-p-1;-} $ is $2p$-Hamilton-biconnected, a contradiction.
Next suppose that $m=n-1$. Since  $ M_{n,n-1}^{n-k,k-p-1;-} $ is a spanning subgraph of  $ M_{n,n-1}^{n-k,k-p-1} $, it follows from Lemma~\ref{R3-Lem2.3}~($i$) that  $ M_{n,n-1}^{n-k,k-p-1;-} $ is not $2p$-Hamilton-biconnected, as desired.

\medskip
{\bf{Case~2}.} $t=m-k+p+1$. Then  $|Y_1|=m-k+p+1$ and $k\geq p+2$.
Set $X=\bigcup _{i=1}^2 X_i$, where $X_1=\{x\in X:d_G(x)=m-k+p+1\}$,  $X_2=\{x\in X: d_G(x)\geq m-k+p+2\}$.
Set $Y_2=Y\backslash Y_1$.
Obviously, $|Y_2|=k-p-1$.  Since  $G$ is an $(m+p+2)$-closed  bipartite graph with $\delta(G)\geq k$, every vertex in $Y_2$ is adjacent to every vertex in $X_2$.
We claim that $Y_2=\{y\in Y:d_G(y)=k\}$. Otherwise, there exists a vertex in $Y_2$ adjacent to every vertex in $X$, and thus $t\geq m-k+p+2$, a contradiction. It follows that $|X_2|=k$ and $|X_1|=n-k$.
Hence $G= M_{n,m}^{n-k,k-p-1} $ for $k\geq p+2$.
Suppose that $m=n$.   Since $M_{n,n}^{n-k,k-p-1;-} $ is a spanning subgraph of $M_{n,n}^{n-k,k-p-1} $,  it follows from  Lemma~\ref{R2-Lem2.3}~($i$) that
 $M_{n,n}^{n-k,k-p-1} $ is   $2p$-Hamilton-biconnected, a contradiction.
Next suppose that $m=n-1$. By Lemma~\ref{R3-Lem2.3}~($i$), $ M_{n,n-1}^{n-k,k-p-1} $ is  not $2p$-Hamilton-biconnected, as desired.
\end{Proof}


\bigskip

Now we are ready to prove Theorems~\ref{Thm1.7} and \ref{Thm1.8}.

\medskip

\noindent {\bf Proof of Theorem~\ref{Thm1.7}.}
Suppose that $ k\geq p$, $e(G)> n(n-k+p-1)+(k+2)(k-p+1)$, and $G$ is not $2p$ Hamilton-biconnected. Let $H=cl_{n+p+2}(G)$.
By Lemma~\ref{Lem3.1}, $H$ is also not $2p$-Hamilton-biconnected. Furthermore, $\delta(H)\geq \delta(G)\geq k$ and $e(H)\geq e(G)> n(n-k+p-1)+(k+2)(k-p+1)$.
By Lemma~\ref{Lem2.4}, $K_{n,n-k+p}\subseteq H$, or $H\in\{ N_{n,n}^{p,1}, N_{n,n}^{p,2}\}$ for $k=p+2$.
It follows from Lemmas~ \ref{R2-Lem2.3}~($ii$), \ref{R3-Lem2.3}~($iii$), and  \ref{Lem2.6}~($i$) that $H=M_{n,n}^{n-k,k-p}  $ for  $k\geq p+1$,
 or $H=  N_{n,n}^{p,1}$ for $k=p+2$. Hence  $G\subseteq M_{n,n}^{n-k,k-p}  $ for  $k\geq p+1$,
 or $G\subseteq  N_{n,n}^{p,1}$ for $k=p+2$. \QEDB
\medskip

Let $p=0$ in Theorem~\ref{Thm1.7}, we partially prove the following Moon and Moser's Theorem \cite{MM}.

\begin{Corollary}\label{Cor4}
Let $G$ be a  balanced bipartite graph of order $2n$ with  $\delta(G)\geq k$, where $1\leq k\leq\frac{n-2}{3}$. If
$$e(G)> n(n-k)+k^2,$$ then $G$ is Hamiltonian.
\end{Corollary}

\begin{Proof}
Note that $e(G)>n(n-k)+k^2\geq n(n-k-1)+(k+2)(k+1)$, $e(M_{n,n}^{n-k,k})=n(n-k)+k^2$, and $e(N_{n,n}^{0,1})=n^2-2n+4$. It follows from Theorem~\ref{Thm1.7} that $G$ is Hamilton-biconnected. Hence $G$ is  Hamiltonian.
\end{Proof}

\bigskip

\noindent {\bf Proof of Theorem~\ref{Thm1.8}.}
Denote $G=(X,Y;E) $  with $|X|=n$ and  $|Y|=n-1$.
Suppose that $ k\geq p$, $e(G)>n(n-k+p-2)+(k+1)(k-p+1)$, and $G$ is not $2p$-Hamilton-biconnected. Let $H=cl_{n+p+1}(G)$.
By Lemma~\ref{Lem4.2}, $H$ is also not $2p$-Hamilton-biconnected. In addition, $\delta(H)\geq \delta(G)\geq k$ and $e(H)\geq e(G)> n(n-k+p-2)+(k+2)(k-p+1)$.
By Corollary~\ref{Cor2.5},   $K_{n,n-k+p-1}\subseteq H$ or  $K_{n-1,n-k+p}\subseteq H$.
Since $H$ is not $2p$-Hamilton-biconnected, we have $H\neq K_{n,n-1}$, which implies that $k \geq p+1$.
 Next we consider the following two cases.

\medskip
{\bf Case~1.} $K_{n,n-k+p-1}\subseteq H$.  Note that $G\subseteq H$, $M_{n,n-1}^{n-k,k-p}\subseteq M_{n,n-1}^{n-k-1,k-p}$, and $M_{n,n-1}^{n-k,k-p-1;-} \subseteq M_{n,n-1}^{n-k,k-p-1}$.
Combining this with Lemma~\ref{Lem2.6}~($ii$), $G$ is $2p$-Hamilton-biconnected unless $G\subseteq M_{n,n-1}^{n-k-1,k-p}$ for $k\geq p+1$, or $G\subseteq M_{n,n-1}^{n-k,k-p-1}$ for $k\geq p+2$.

\medskip
{\bf Case~2.} $K_{n-1,n-k+p}\subseteq H$ and $K_{n,n-k+p-1}\nsubseteq H$.
Let  $s$, $t$ with $s\geq t$ be the largest integers such that $K_{s,t}\subseteq H$.  It follows that
$s=n-1$ and $n-k+p\leq t\leq n-1$.   We consider the following two subcases.

\medskip
 {\bf Case~2.1.}  Let $X_1\subseteq X$ with $|X_1|=t$ such that $H[X_1, Y]=K_{t,n-1}$.  We show that $ t=n-k+p$. 
Indeed, if $ t>n-k+p$, then $d_H(y)\geq n-k+p+1$ for every $y\in Y$. Since $H$ is an $(n+p+1)$-closed bipartite graph with $\delta(H)\geq k$,
every vertex in $Y$ is adjacent to every vertex in $X$  and thus $H=K_{n,n-1}$, a contradiction.
Then $|X_1|=n-k+p$ and $|X\backslash X_1|=k-p$.
Furthermore, since $H$ is an $(n+p+1)$-closed bipartite graph with $\delta(H)\geq k$, $d_H(x)=k$ for every $x\in X\backslash X_1$.
Let $Y_1=\{y\in Y:d_H(y)=n-k+p\}$. Moreover,   every vertex in $ X\backslash X_1$ is adjacent to every vertex in $Y\backslash Y_1$.
It follows that $|Y\backslash Y_1|=k$ and  $|Y_1|=n-k-1$.
Hence $H= M_{n,n-1}^{k-p,n-k-1}$ for $k\geq p+1$.   On the other hand, by Lemma~\ref{R3-Lem2.3}~($i$), $M_{n,n-1}^{k-p,n-k-1}$ is  not $2p$-Hamilton-biconnected, as desired.

\medskip
 {\bf Case~2.2.}  Let $X_1\subseteq X$ and $Y_1\subseteq Y$ with $|X_1|=n-1$ and $|Y_1|=t$ such that $H[X_1, Y_1]=K_{n-1,t}$.
We first show that $k=p+1$.
Since $H$ is an  $(n+p+1)$-closed bipartite graph, if $k>p+1$ then
  every vertex in $X$ is adjacent to every vertex in  $Y_1$. This implies
that $K_{n,n-k+p-1}\subseteq K_{n,n-k+p}\subseteq H$, a contradiction.
Since $k=p+1$, we have $t=n-1$. Hence $K_{n-1,n-1}\subseteq H$, which can be described to Case~2.1. \QEDB

\section{Proofs of Theorems~\ref{Thm1.9} and \ref{Thm1.10}}
In order to prove Theorems~\ref{Thm1.9} and \ref{Thm1.10}, we need the following lemma.

\begin{Lemma}\label{Lem3.2}
(i) For  $k\geq2$ and $n\geq 2k^2+3$, $\rho(N_{n,n}^{k-2,1})> \rho(M_{n,n}^{n-k,2})$.\\
(ii) For  $k\geq2$ and $n\geq k+1$, $q(M_{n,n}^{n-k,2})>q(N_{n,n}^{k-2,1})$.
\end{Lemma}

\begin{Proof}
(i) Denote $M_{n,n}^{n-k,2}=(X,Y;E)$  with $|X|=|Y|=n$.
Let $\mathbf{x}$ be the eigenvector corresponding to  $\rho(M_{n,n}^{n-k,2})$.
 Let $X=X_1\bigcup X_2$ and $Y=Y_1\bigcup Y_2$, where $X_1$ and $X_2$ are the sets of vertices in $X$ with degree $n-2$ and $n$ respectively, and  $Y_1$ and $Y_2$ are the sets of vertices in $Y$ with degree $n$ and $k$ respectively.

By symmetry, the entry of $\mathbf{x}$ corresponding to any vertex in $X_i$, denoted by $x_i$,  is a positive constant for $1\leq i\leq2$.
Similarly, the entry of $\mathbf{x}$ corresponding to any vertex in $Y_i$, denoted by $y_i$, is also a positive constant for $1\leq i\leq2$.
By eigenequation $A(M_{n,n}^{n-k,2})\mathbf{x}=\rho( M_{n,n}^{n-k,2})\mathbf{x}$, we have
\begin{eqnarray*}
  \rho x_1 &=& (n-2)y_1, \\
  \rho x_2 &=& (n-2)y_1+2y_2, \\
  \rho y_1 &=& (n-k)x_1+kx_2,\\
  \rho y_2 &=& kx_2.
\end{eqnarray*}
By a simple calculation, $\rho(M_{n,n}^{n-k,2})$ is the largest root of $f(x)=0$, where
 $$f(x)=x^4-(n^2-2n+2k)x^2+2k(n-k)(n-2).$$
Since $$f(n-1)=n(n-2k^2-2)+4k^2-2k+1>0$$
 and for $x\geq n-1$,
 $$f'(x)=2x(2x^2-n^2+2n-2k)\geq 2(n-1)(n^2-2n-2k+2)>0,$$  we have
$\rho(M_{n,n}^{n-k,2})<n-1$.
On the other hand, since $K_{n-1,n-1}$ is a subgraph of $N_{n,n}^{k-2,1}$, it follows from Lemma~\ref{Lem2.1} that
$$ \rho(N_{n,n}^{k-2,1})\geq \rho(K_{n-1,n-1})=n-1>\rho(M_{n,n}^{n-k,1}).$$

\medskip
(ii) Let $f(x)=x(x-n)f_1(x)$ and $g(x)=x(x-n)^2(x-k+1)g_1(x)$,
where $$f_1(x)=x^2-(2n+k-2)x+2kn-4k,$$ $$ g_1(x)=x^2-(2n+k-1)x+2kn+2n-4k.$$
By a similar argument to the proof of ($i$), $q(M_{n,n}^{n-k,2})$ and $q(N_{n,n}^{k-2,1})$
  are the largest roots of $f(x)=0$ and  $g(x)=0$,  respectively.
 Furthermore,  since $K_{n,n-2}$ and $K_{n-1,n-1}$ are  proper subgraphs of $M_{n,n}^{n-k,2}$ and $N_{n,n}^{k-2,1}$, respectively, it follows from Lemma~\ref{Lem2.1} that
 $$q(M_{n,n}^{n-k,2})>2n-2, \quad q(M_{n,n}^{n-k,2})>2n-2.$$ Hence $q(M_{n,n}^{n-k,2})$ and $q(N_{n,n}^{k-2,1})$
  are the largest roots of $f_1(x)=0$ and  $g_1(x)=0$,  respectively.
On the other hand,  since both $M_{n,n}^{n-k,2}$ and $N_{n,n}^{k-2,1}$ are  proper subgraphs of $K_{n,n}$, it follows from Lemma~\ref{Lem2.1} that
$$q(M_{n,n}^{n-k,2})<2n, \quad q(M_{n,n}^{n-k,2})<2n.$$
Since  for  $x<2n$
$$g_1(x)-f_1(x)=2n-x>0,$$ we have $q(M_{n,n}^{n-k,2})>q(N_{n,n}^{k-2,1})$.
\end{Proof}

\bigskip
%
\noindent {\bf Proof of Theorem~\ref{Thm1.9}.}
(i) Suppose that $\rho(G)\geq \rho(N_{n,n}^{k-2,1})$ and $G$ is not $2p$-Hamilton-biconnected.
Since $K_{n-1,n-1} $ is a proper subgraph of $ N_{n,n}^{k-2,1}$, Lemma~\ref{Lem2.1} implies that
$$\rho(G)\geq \rho(N_{n,n}^{k-2,1})>\rho (K_{n-1,n-1})=n-1.$$
By Lemma~\ref{Lem2.2}, $ \sqrt{e(G)}\geq  \rho(G)> n-1$, which implies that
$$e(G)> (n-1)^2\geq n(n-3)+3(k+2).$$
 It follows from Theorem~\ref{Thm1.7} that $G\subseteq M_{n,n}^{n-k,2} $ or $G\subseteq N_{n,n}^{k-2,1}$.
If $G\subseteq M_{n,n}^{n-k,2} $, then Lemmas~\ref{Lem2.1} and ~\ref{Lem3.2}~($i$)
imply that  $\rho(G)\leq\rho(M_{n,n}^{n-k,2})<\rho (  N_{n,n}^{k-2,1}) $, a contradiction.
If $G$ is a proper subgraph of $N_{n,n}^{k-2,1} $, then
 Lemma~\ref{Lem2.1} implies that $\rho(G)<\rho (N_{n,n}^{k-2,1}) $, a contradiction.
Hence $G= N_{n,n}^{k-2,1}$.

\medskip
(ii) Suppose that $\rho(G)\geq \rho(M_{n,n}^{n-k,k-p})$ and  $G$ is not $2p$-Hamilton-biconnected.
 Since $K_{n,n-k+p} $ is a proper subgraph of $ M_{n,n}^{n-k,k-p}$,  Lemma~\ref{Lem2.1} implies that
$$\rho(G)\geq \rho(M_{n,n}^{n-k,k-p})>\rho (K_{n,n-k+p})=\sqrt{n(n-k+p)}.$$
By Lemma~\ref{Lem2.2}, $ \sqrt{e(G)}\geq  \rho(G)> \sqrt{n(n-k+p)}$, which implies that
$$e(G)>n(n-k+p)\geq n(n-k+p-1)+(k+2)(k-p+1).$$
It follows from Theorem~\ref{Thm1.7} that $G\subseteq M_{n,n}^{n-k,k-p} $.
If $G$ is a proper subgraph of $M_{n,n}^{n-k,k-p}$, then Lemma~\ref{Lem2.1} implies  that  $\rho(G)<\rho(M_{n,n}^{n-k,k-p})$, a contradiction.
Hence $G=  M_{n,n}^{n-k,k-p}$. \QEDB

\begin{Corollary}
Let  $p\geq0$ and  $G$ be a  balanced bipartite graph of order $2n$  with $\delta(G)\geq k$, where  $n\geq n_0(k,p)$ and
\[n_0(k,p)=\left\{
\begin{array}{cccc}
 \vspace{1mm}
   2k^2+3,&& \mbox{if $k=p+2$}\\
   \vspace{1mm}
   (k+2)(k-p+1),&& \mbox{otherwise.}
\end{array}\right.
\]
If $k\geq p+1$ and $\rho(G)\geq\sqrt{n(n-k+p)+k(k-p)}$,  then $G$ is  $2p$-Hamilton-biconnected.
\end{Corollary}

\begin{Proof}
Suppose that $k=p+2$. Note that $e(N_{n,n}^{k-2,1})=n^2-2n+2k$. By Lemma~\ref{Lem2.2} and Theorem~\ref{Thm1.9}~($i$), the result follows.
Next suppose that  $k\neq p+2$.  Note that $e(M_{n,n}^{n-k,k-p})=n(n-k+p)+k(k-p)$. By Lemma~\ref{Lem2.2} and Theorem~\ref{Thm1.9}~($ii$), the result follows.
\end{Proof}

\bigskip


\noindent {\bf Proof of Theorem~\ref{Thm1.10}.}
Suppose that $q(G)\geq q(M_{n,n}^{n-k,k-p})$ and  $G$ is not $2p$-Hamilton-biconnected.
 Since $K_{n,n-k+p}$ is a proper subgraph of $ M_{n,n}^{n-k,k-p}$, Lemma~\ref{Lem2.1} implies that
$$q(G)\geq q( M_{n,n}^{n-k,k-p})> q (K_{n,n-k+p})=2n-k+p.$$
By Lemma~\ref{Lem2.3}, $\frac{e(G)}{n}+n\geq q(G)> 2n-k+p $,  which implies that
$$e(G)> n(n-k+p)\geq n(n-k+p-1)+(k+2)(k-p+1).$$
It follows from Theorem~\ref{Thm1.7} that $G\subseteq M_{n,n}^{n-k,k-p}  $ for  $k\geq p+1$, or $G\subseteq  N_{n,n}^{p,1}$ for $k=p+2$.
If $G$ is a proper subgraph of $M_{n,n}^{n-k,k-p}$, then Lemma~\ref{Lem2.1} implies that $q(G)<q (M_{n,n}^{n-k,k-p})$, a contradiction.
If $G$ is a subgraph of $N_{n,n}^{p,1}$ for $k=p+2$, then Lemmas~\ref{Lem2.1} and \ref{Lem3.2}~($ii$) imply that
$q(G)\leq q (N_{n,n}^{p,1})<q (M_{n,n}^{n-k,k-p})$,  a contradiction. Hence $G=  M_{n,n}^{n-k,k-p}$.
\QEDB

\begin{Corollary}
Let  $p\geq0$ and $G$ be a  balanced bipartite graph of order $2n$ with  $\delta(G)\geq k$, where $n\geq (k+2)(k-p+1)$.
If  $k\geq p+1$ and  $q(G)\geq 2n-k+p+\frac{k(k-p)}{n}$,  then $G$ is $2p$-Hamilton-biconnected.
\end{Corollary}

\begin{Proof}
Note that $n+\frac{e(M_{n,n}^{n-k,k-p})}{n}=2n-k+p+\frac{k(k-p)}{n}$. By Lemma~\ref{Lem2.3} and Theorem~\ref{Thm1.10},  the result follows.
\end{Proof}

\section{Proofs of Theorems~\ref{Thm1.11} and \ref{Thm1.12}}
The proofs of Lemmas~\ref{Lem4.3} and \ref{Lem4.4} are similar to that of  Lemma~\ref{Lem3.2}, so we put them in the  appendix.

\begin{Lemma}\label{Lem4.3}
 (i) For  $p\geq0$, $k\geq p+1$,  and  $n\geq 2k-p+2$,
 $$\rho(M_{n,n-1}^{k-p,n-k-1})>\rho(M_{n,n-1}^{n-k-1,k-p}).$$

(ii) For  $p\geq0$,  $k\geq p+2$, and $n\geq 2k-p+2$,
$$\rho(M_{n,n-1}^{n-k,k-p-1})> \rho(M_{n,n-1}^{k-p,n-k-1}).$$
\end{Lemma}

\begin{Lemma}\label{Lem4.4}
(i) For $p\geq0$,  $k\geq p+1$,  and $n\geq 2k-p+2$,
 $$q(M_{n,n-1}^{n-k-1,k-p})>q(M_{n,n-1}^{k-p,n-k-1}).$$

 (ii) For  $p\geq0$, $k\geq p+2$, and  $n\geq 2k-p+2$,
$$q(M_{n,n-1}^{n-k,k-p-1})> q(M_{n,n-1}^{n-k-1,k-p}).$$
\end{Lemma}



\noindent {\bf Proof of Theorem~\ref{Thm1.11}.}
(i) Suppose that $\rho(G)\geq \rho(M_{n,n-1}^{1,n-k-1})$ and  $G$ is not $2p$-Hamilton-biconnected.
 Since $K_{n-1,n-1}$ is a proper subgraph of $M_{n,n-1}^{1,n-k-1}$,
it follows from Lemma~\ref{Lem2.1} that
$$ \rho(G)\geq\rho(M_{n,n-1}^{1,n-k-1})>\rho(K_{n-1,n-1})=n-1.$$
By  Lemma~\ref{Lem2.2} $\sqrt{e(G)}\geq  \rho(G)>n-1 $, which implies that
$$e(G)> n^2-2n+1\geq n(n-3)+2k+4.$$
Then it follows from  Theorem~\ref{Thm1.8}  that   $G\subseteq M_{n,n-1}^{n-k-1,1} $ or $G\subseteq M_{n,n-1}^{1,n-k-1}$.
By Lemmas~\ref{Lem2.1} and \ref{Lem4.3}~($i$), $G= M_{n,n-1}^{1,n-k-1}$.

\medskip
(ii) Suppose that  $\rho(G)\geq \rho(M_{n,n-1}^{n-k,k-p-1})$ and  $G$ is not $2p$-Hamilton-biconnected.
 Since $K_{n,n-k+p}$ is a proper subgraph of $M_{n,n-1}^{n-k,k-p-1}$, it follows from Lemma~\ref{Lem2.1} that
$$\rho(G)\geq \rho(M_{n,n-1}^{n-k,k-p-1})>\rho (K_{n,n-k+p})=\sqrt{n(n-k+p)}.$$
By Lemma~\ref{Lem2.2}, $ \sqrt{e(G)}\geq  \rho(G)>\sqrt{n(n-k+p)}$, which implies that
$$e(G)> n(n-k+p)\geq n(n-k+p-2)+(k+2)(k-p+1).$$
It follows from  Theorem~\ref{Thm1.8}  that  $G\subseteq M_{n,n-1}^{n-k-1,k-p} $,  $G\subseteq M_{n,n-1}^{k-p,n-k-1}$, or
 $G\subseteq M_{n,n-1}^{n-k,k-p-1}$.
By Lemmas~\ref{Lem2.1} and \ref{Lem4.3}, $G= M_{n,n-1}^{n-k,k-p-1}$. \QEDB

\begin{Corollary}
Let  $p\geq0$ and $G$ be a nearly balanced bipartite graph of order $2n-1$ and $\delta(G)\geq k$.  \\
(i) $k=p+1$, $n\geq2k+3$, and  $\rho(G)\geq\sqrt{(n-1)^2+k}$, then $G$ is  $2p$-Hamilton-biconnected.\\
(ii) If $k\geq p+2$, $n\geq \frac{(k+2)(k-p+1)}{2}$,  and $\rho(G)\geq\sqrt{n(n-k+p)+k(k-p-1)}$,   then $G$ is $2p$-Hamilton-biconnected.
\end{Corollary}

\begin{Proof}
(i) Note that $e(M_{n,n-1}^{1,n-k-1})=(n-1)^2+k$. By Lemma~\ref{Lem2.2} and Theorem~\ref{Thm1.11}~($i$), the result follows.

\medskip
(ii) Note that $e(M_{n,n-1}^{n-k,k-p-1})=n(n-k+p)+k(k-p-1)$. By Lemma~\ref{Lem2.2} and Theorem~\ref{Thm1.11}~($ii$), the result follows.
\end{Proof}

\bigskip


\noindent {\bf Proof of Theorem~\ref{Thm1.12}.}
(i) Suppose that $q(G)\geq q(M_{n,n-1}^{n-k-1,1})$ and  $G$ is not $2p$-Hamilton-biconnected.
 Since $K_{n,n-2}$ is a proper subgraph of $M_{n,n-1}^{n-k-1,1}$,
 Lemma~\ref{Lem2.1} implies that
$$ q(G)\geq q(M_{n,n-1}^{n-k-1,1})>q(K_{n,n-2})=2n-2.$$
By  Lemma~\ref{Lem2.3}, $\frac{ e(G)}{n}+n\geq  q(G)>2n-2$.
Note that here we consider $G$ as a balanced bipartite graph having an isolated vertex. This implies that
$$e(G)> n^2-2n\geq n(n-3)+2k+4.$$
It follows from  Theorem~\ref{Thm1.8}  that   $G\subseteq M_{n,n-1}^{n-k-1,1} $ or $G\subseteq M_{n,n-1}^{1,n-k-1}$.
By Lemmas~\ref{Lem2.1} and \ref{Lem4.4}~($i$), $G= M_{n,n-1}^{1,n-k-1}$.

\medskip
(ii) Suppose that  $q(G)\geq q(M_{n,n-1}^{n-k,k-p-1})$ and  $G$ is not $2p$-Hamilton-biconnected.
  Since $K_{n,n-k+p}$  is a proper subgraph of  $M_{n,n-1}^{n-k,k-p-1}$,  Lemma~\ref{Lem2.1} implies that
$$q(G)\geq q(M_{n,n-1}^{n-k,k-p-1})> q (K_{n,n-k+p})=2n-k+p.$$
By Lemma~\ref{Lem2.3}, $\frac{e(G)}{n}+n\geq q(G)>2n-k+p $. Note that here we consider $G$ as a balanced bipartite graph having an isolated vertex.
This implies that
$$e(G)>n(n-k+p)\geq n(n-k+p-2)+(k+2)(k-p+1).$$
It follows from  Theorem~\ref{Thm1.8}  that  $G\subseteq M_{n,n-1}^{n-k-1,k-p} $, $G\subseteq M_{n,n-1}^{n-k,k-p-1}$,  or $G\subseteq M_{n,n-1}^{k-p,n-k-1}$.
By Lemmas~\ref{Lem2.1} and \ref{Lem4.4}, $G= M_{n,n-1}^{n-k,k-p-1}$. \QEDB

\begin{Corollary}
Let  $p\geq0$ and $G$ be a nearly balanced bipartite graph of order $2n-1$ and $\delta(G)\geq k$.  \\
(i) If  $k=p+1$,  $n\geq2k+4$,and $q(G)\geq \sqrt{2n-2+\frac{k+1}{n}}$,   then $G$ is  $2p$-Hamilton-biconnected.\\
(ii) If $k\geq p+2$, $n\geq \frac{(k+2)(k-p+1)}{2}$, and  $q(G)\geq \sqrt{2n-k+p+\frac{k(k-p-1)}{n}}$, then $G$ is $2p$-Hamilton-biconnected.
\end{Corollary}

\begin{Proof}
(i) Note that $n+\frac{e(M_{n,n-1}^{n-k-1,1})}{n}=2n-2+\frac{k+1}{n}$. By Lemma~\ref{Lem2.3} and Theorem~\ref{Thm1.12}~($i$), the result follows.

\medskip
(ii) Note that $n+\frac{e(M_{n,n-1}^{n-k,k-p-1})}{n}=2n-k+p+\frac{k(k-p-1)}{n}$. By Lemma~\ref{Lem2.3} and Theorem~\ref{Thm1.12}~($ii$), the result follows.
\end{Proof}

\vspace{16mm}
\newpage

\noindent {\bf  \large Appendix}

\medskip

Denote by $P_{uv}$  a path between $u$ and $v$.
Denote  by   $P_{uv} \bigsqcup P_{wz}$ a path obtained from two disjoint paths $P_{uv}$ and $P_{wz}$ by joining $v$ and $w$.

\bigskip

\begin{tikzpicture}[scale=1.8]
\path (1.5,0)  coordinate (P1);\path (2.2,0)  coordinate (P2);  \path (2.8,0)  coordinate (P3); \path (3.5,0)  coordinate (P4);
\path (1.5,-1)  coordinate (Q1);\path (2.2,-1)  coordinate (Q2);  \path (2.8,-1)  coordinate (Q3); \path (3.5,-1)  coordinate (Q4);
\path (2.5,0)  coordinate (S1);\path (2.5,-1)  coordinate (S2);
\path (1.7,0)  coordinate (R1);\path (1.85,0)  coordinate (R2);  \path (2,0)  coordinate (R3);
\path (3,0)  coordinate (R4); \path (3.15,0)  coordinate (R5);\path (3.3,0)  coordinate (R6);
\path (1.7,-1)  coordinate (R7);\path (1.85,-1)  coordinate (R8);  \path (2,-1)  coordinate (R9);
\path (3,-1)  coordinate (R10); \path (3.15,-1)  coordinate (R11);\path (3.3,-1)  coordinate (R12);
\fill (S1) circle (1.5pt);\fill (S2) circle (1.5pt);
 \foreach \i in {1,2,3,4}
{\fill (P\i) circle (1.5pt);\fill (Q\i) circle (1.5pt);}
\foreach \i in {1,...,12}
{\fill (R\i) circle (0.5pt);}
 \foreach \i in {1,2}
 { \draw (S1)--(Q\i);\foreach \j in {1,2,3,4}
{\draw (Q\i)--(P\j);}}
 \foreach \i in {3,4}
 {\draw (S2)--(P\i); \foreach \j in {3,4}
{\draw (P\i)--(Q\j);}}
\draw (S1)--(S2);
\path (1.5,0.2) node () {$u_{11}$};\path (2,0.2) node () {$u_{1,s-1}$};
\path (2.5,0.2) node () {$u$};\path (2.5,-1.2) node () {$v$};
\path (2.8,0.2) node () {$u_{21}$};\path (3.4,0.2) node () {$u_{2,t+1}$};
\path (1.5,-1.2) node () {$v_{11}$};\path (2.2,-1.2) node () {$v_{1s}$};
\path (2.8,-1.2) node () {$v_{21}$};\path (3.5,-1.2) node () {$v_{2t}$};
\path (2.5,-1.6) node () {$M_{s+t+1,s+t+1}^{s,t;-}$};

\path (4,0)  coordinate (P1);\path (4.5,0)  coordinate (P2); \path (4.9,0)  coordinate (P3);\path (5.4,0)  coordinate (P4);
 \path (6,0)  coordinate (P5); \path (6.5,0)  coordinate (P6); \path (6.7,0)  coordinate (P7); \path (7.2,0)  coordinate (P8);
\path (4.4,-1)  coordinate (Q1);\path (5.2,-1)  coordinate (Q2);  \path (5.9,-1)  coordinate (Q3); \path (6.7,-1)  coordinate (Q4);
\path (4.15,0)  coordinate (R1);\path (4.25,0)  coordinate (R2);  \path (4.35,0)  coordinate (R3);
\path (5.05,0)  coordinate (R4); \path (5.15,0)  coordinate (R5);\path (5.25,0)  coordinate (R6);
\path (6.15,0)  coordinate (R7);\path (6.25,0)  coordinate (R8);  \path (6.35,0)  coordinate (R9);
\path (6.85,0)  coordinate (R10); \path (6.95,0)  coordinate (R11);\path (7.05,0)  coordinate (R12);
\path (4.6,-1)  coordinate (R13); \path (4.8,-1)  coordinate (R14);\path (5,-1)  coordinate (R15);
\path (6.1,-1)  coordinate (R16); \path (6.3,-1)  coordinate (R17);\path (6.5,-1)  coordinate (R18);
\fill (S1) circle (1.5pt);\fill (S2) circle (1.5pt);
 \foreach \i in {1,...,8}
{\fill (P\i) circle (1.5pt);}
 \foreach \i in {1,...,4}
{\fill (Q\i) circle (1.5pt);}
\foreach \i in {1,...,18}
{\fill (R\i) circle (0.5pt);}
 \foreach \i in {1,2}
 {\foreach \j in {1,...,8}
{\draw (Q\i)--(P\j);}}
 \foreach \i in {3,4}
 { \foreach \j in {3,4}
{\draw (P\i)--(Q\j);}}
\draw (Q3)--(P5);\draw (Q3)--(P6);\draw (Q4)--(P7);\draw (Q4)--(P8);
\path (4,0.2) node () {$u_{11}$};\path (4.4,0.2) node () {$u_{1t}$};
\path (4.8,0.2) node () {$u_{21}$};\path (5.3,0.2) node () {$u_{2g}$};
\path (5.8,0.3) node () {$u_{31}^{(1)}$};\path (6.3,0.3) node () {$u_{3,s+1}^{(1)}$};
\path (6.8,0.4) node () {$u_{31}^{(h)}$};\path (7.4,0.4) node () {$u_{3,s+1}^{(h)}$};
\path (4.4,-1.2) node () {$v_{11}$};\path (5.2,-1.2) node () {$v_{1,m-k}$};
\path (5.9,-1.2) node () {$v_{21}$};\path (6.7,-1.2) node () {$v_{2h}$};
\path (5.5,-1.6) node () {$H_{n-p,m-p}^{k,p,l}$};
\end{tikzpicture}

\begin{tikzpicture}[scale=1.8]
\path (1.5,0)  coordinate (P1);\path (2.2,0)  coordinate (P2);  \path (2.6,0)  coordinate (P3);\path (3,0)  coordinate (P4);
 \path (3.5,0)  coordinate (P5);\path (3.5,-1)  coordinate (P6);
\path (1.5,-1)  coordinate (Q1);\path (2.2,-1)  coordinate (Q2);  \path (2.6,-1)  coordinate (Q3); \path (3,-1)  coordinate (Q4);
\path (2.5,0)  coordinate (S1);\path (2.5,-1)  coordinate (S2);
\path (1.7,0)  coordinate (R1);\path (1.85,0)  coordinate (R2);  \path (2,0)  coordinate (R3);
\path (1.7,-1)  coordinate (R4);\path (1.85,-1)  coordinate (R5);  \path (2,-1)  coordinate (R6);
\fill (P5) circle (1.5pt);\fill (P6) circle (1.5pt);
 \foreach \i in {1,2,3,4}
{\fill (P\i) circle (1.5pt);\fill (Q\i) circle (1.5pt);}
\foreach \i in {1,...,6}
{\fill (R\i) circle (0.5pt);}
 \foreach \i in {1,2,3,4}
 { \foreach \j in {1,2,3,4}
{\draw (Q\i)--(P\j);}}
 \foreach \i in {3,4}
 { \foreach \j in {3,4}
{\draw (P\i)--(Q\j);}}
\draw (P5)--(Q3);\draw (P5)--(Q4);\draw (P6)--(P3);\draw (P6)--(P4);
\path (1.5,0.2) node () {$u_{11}$};\path (2,0.2) node () {$u_{1t}$};
\path (2.6,0.2) node () {$u_{21}$};\path (3,0.2) node () {$u_{22}$};
\path (1.5,-1.2) node () {$v_{11}$};\path (2.2,-1.2) node () {$v_{1t}$};
\path (2.6,-1.2) node () {$v_{21}$};\path (3,-1.2) node () {$v_{22}$};
\path (3.5,0.2) node () {$u$};\path (3.5,-1.2) node () {$v$};
\path (2.5,-1.6) node () {$N_{n-p,n-p}^{0,2}$};

 \path (4.3,0)  coordinate (P1);\path (5.2,0)  coordinate (P2);
 \path (5.8,0)  coordinate (P3); \path (6.3,0)  coordinate (P4); \path (6.7,0)  coordinate (P5); \path (7.2,0)  coordinate (P6);
\path (4.4,-1)  coordinate (Q1);\path (5.2,-1)  coordinate (Q2);  \path (5.9,-1)  coordinate (Q3); \path (6.7,-1)  coordinate (Q4);
\path (4.5,0)  coordinate (R1); \path (4.7,0)  coordinate (R2);\path (4.9,0)  coordinate (R3);
\path (5.95,0)  coordinate (R4);\path (6.05,0)  coordinate (R5);  \path (6.15,0)  coordinate (R6);
\path (6.85,0)  coordinate (R7); \path (6.95,0)  coordinate (R8);\path (7.05,0)  coordinate (R9);
\path (4.6,-1)  coordinate (R10); \path (4.8,-1)  coordinate (R11);\path (5,-1)  coordinate (R12);
\path (6.1,-1)  coordinate (R13); \path (6.3,-1)  coordinate (R14);\path (6.5,-1)  coordinate (R15);
 \foreach \i in {1,...,6}
{\fill (P\i) circle (1.5pt);}
 \foreach \i in {1,...,4}
{\fill (Q\i) circle (1.5pt);}
\foreach \i in {1,...,18}
{\fill (R\i) circle (0.5pt);}
 \foreach \i in {1,2}
 {\foreach \j in {1,...,8}
{\draw (Q\i)--(P\j);}}
 \foreach \i in {3,4}
 { \foreach \j in {3,4}
{\draw (P\i)--(Q\j);}}
\draw (Q3)--(P5);\draw (Q3)--(P6);\draw (Q4)--(P7);\draw (Q4)--(P8);
\path (4.3,0.2) node () {$u_{11}$};\path (5.2,0.2) node () {$u_{1t}$};
\path (5.8,0.3) node () {$u_{31}^{(1)}$};\path (6.3,0.3) node () {$u_{3,s_1+1}^{(1)}$};
\path (6.8,0.4) node () {$u_{31}^{(h)}$};\path (7.4,0.4) node () {$u_{3,s_h+1}^{(h)}$};
\path (4.4,-1.2) node () {$v_{11}$};\path (5.2,-1.2) node () {$v_{1,m-k}$};
\path (5.9,-1.2) node () {$v_{21}$};\path (6.7,-1.2) node () {$v_{2h}$};
\path (5.5,-1.6) node () {$H_{n-p,m-p}^{k,p,l}(s_1,\ldots,s_{h})$};
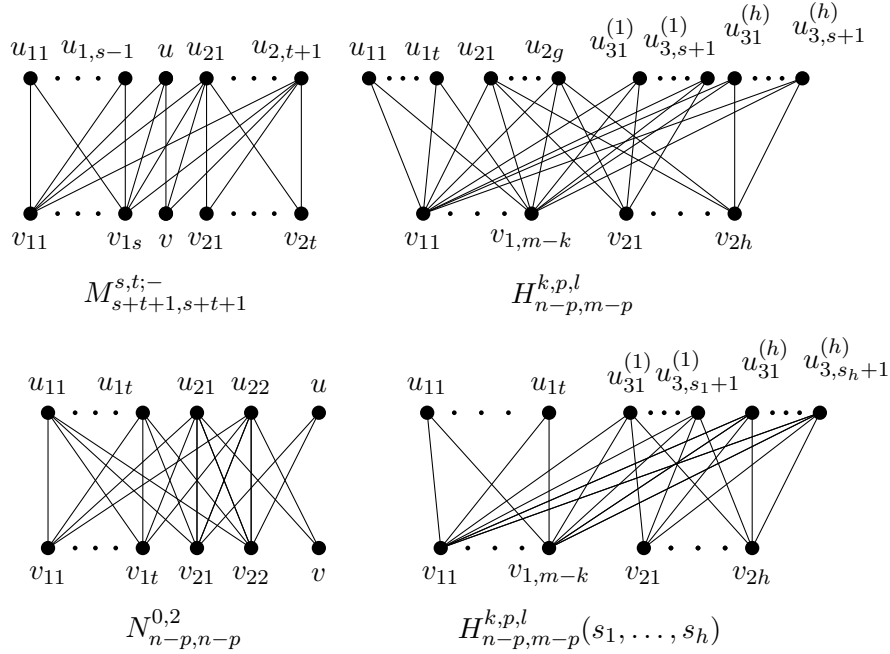
\path (5,-2.2) node () {Fig.~2. Graphs $M_{s+t+1,s+t+1}^{s,t}$, $M_{s+t+1,s+t+1}^{s,t;-}$, $N_{n-p,n-p}^{0,2}$, and $H_{n-p,m-p}^{k,p,l}(s_1,\ldots,s_{h})$.};
\end{tikzpicture}

\bigskip

 \noindent{\bf Proof of Lemma~\ref{R-Lem2.3}.}
We first assume that $l>p$.   Let  $g=l-p$, $h=k-p$, $s=k-l-1$, and $t=n-(k-p)(k-l)-l$. Let  $H_{n-p,m-p}^{k,p,l}$ be a bipartite graph obtained from $F_{n,m}^{k,p,l}$ by deleting all vertices in a balanced set of size $2p$ which consists of vertices with as large as possible degree (see Fig.~2).  Note that every bipartite graph of order $m+n-2p$  obtained from $F_{n,m}^{k,p,l}$ by deleting all vertices in a balanced set of size $2p$ contains $H_{n-p,m-p}^{k,p,l}$ as a subgraph.  It suffices to prove that $H_{n-p,m-p}^{k,p,l}$ is Hamilton-biconnected.
Label the vertices of $H_{n-p,m-p}^{k,p,l}$ as $u_{11},\ldots,u_{1t}$, $u_{21},\ldots,u_{2g}$, $u_{31}^{(1)},\ldots,u_{3,s+1}^{(1)},\ldots, u_{31}^{(h)},\ldots,u_{3,s+1}^{(h)}$, $v_{11},\ldots,v_{1,m-k}$, $v_{21},\ldots,v_{2h}$ (see Fig.~2).

Let $m=n$.  We assume that $l\leq k-2$. Clearly, $g\geq1$, $h\geq2$, and $s\geq1$. Denote
\begin{eqnarray*}
  P_{u_{11}v_{1t}} &=& \bigsqcup_{i=1}^t u_{1i}v_{1i},  \quad P_{u_{21}v_{1,t+g}}=\bigsqcup_{i=1}^{g} u_{2i}v_{1,t+i},\\
 P_{u_{31}^{(h)}v_{2h}} &=& \bigsqcup_{i=1}^{s} u_{3i}^{(h)}v_{1,t+g+(h-1)s+i}\bigsqcup u_{3,s+1}^{(h)}v_{2h}, \\
  Q_{u_{31}^{(h)}v_{2h}} &=& \bigsqcup_{i=1}^{s-1} u_{3i}^{(h)}v_{1,t+g+(h-1)s+i}\bigsqcup u_{3s}^{(h)}v_{2h},\\
    P_{u_{3,s+1}^{(i)}v_{1,t+g+is}}&=& u_{3,s+1}^{(i)}v_{2i}\bigsqcup \bigg( \bigsqcup_{j=1}^{s} u_{3j}^{(i)}v_{1,t+g+(i-1)s+j}\bigg)~ \mbox{for $1\leq i\leq h$}.\\
\end{eqnarray*}
%
%
%
$H_{n-p,n-p}^{k,p,l}$ has seven kinds of Hamiltonian paths, denoted by  $R_1, \ldots, R_7$. We present them as follows:
\begin{eqnarray*}
  R_1 &=& P_{u_{11}v_{1t}}\bigsqcup P_{u_{21}v_{1,t+g}}\bigsqcup \bigg(\bigsqcup_{i=1}^{h} P_{u_{3,s+1}^{(i)}v_{1,t+g+is}}\bigg), \\
 R_2 &=& P_{u_{11}v_{1t}} \bigsqcup P_{u_{21}v_{1,t+g}} \bigsqcup \bigg(\bigsqcup_{i=1}^{h-1} P_{u_{3,s+1}^{(i)}v_{1,t+g+is}} \bigg)\bigsqcup P_{u_{3,1}^{(h)}v_{2,h}}, \\
 R_3 &=&P_{u_{21}v_{1,t+g}} \bigsqcup P_{u_{11}v_{1t}}\bigsqcup \bigg(\bigsqcup_{i=1}^{h} P_{u_{3,s+1}^{(i)}v_{1,t+g+is}}\bigg), \\
  R_4 &=&P_{u_{21}v_{1,t+g}}\bigsqcup P_{u_{11}v_{1t}}\bigsqcup \bigg(\bigsqcup_{i=1}^{h-1} P_{u_{3,s+1}^{(i)}v_{1,t+g+is}} \bigg)\bigsqcup P_{u_{3,1}^{(h)}v_{2h}}, \\
  R_5 &=& \bigsqcup_{i=1}^{h} P_{u_{3,s+1}^{(i)}v_{1,t+g+is}}\bigsqcup P_{u_{11}v_{1t}}\bigsqcup P_{u_{21}v_{1,t+g}}, \\
  R_6 &=& \bigsqcup_{i=1}^{h-1} P_{u_{3,s+1}^{(i)}v_{1,t+g+is}} \bigsqcup P_{u_{11}v_{1t}}\bigsqcup P_{u_{21}v_{1,t+g}} \bigsqcup P_{u_{31}^{(h)}v_{2h}}, \\
  R_7 &=& u_{3,s+1}^{(h)}v_{1,n-k}\bigsqcup \bigg(\bigsqcup_{i=1}^{h-1} P_{u_{3,s+1}^{(i)}v_{1,t+g+is}}\bigg)\bigsqcup P_{u_{11}v_{1t}} \bigsqcup P_{u_{21}v_{1,t+g}}\bigsqcup Q_{u_{31}^{(h)}v_{2h}}.
\end{eqnarray*}
Hence $H_{n-p,n-p}^{k,p,l}$ is Hamilton-biconnected.
Thus $F_{n,n}^{k,p,l}$ is  $2p$-Hamilton-biconnected for $ p< l\leq k-2$.   Similarly we can prove that
$F_{n,n}^{k,p,l}$ is also $2p$-Hamilton-biconnected for $p<l= k-1$.

Let $m=n-1$.    We assume that $l\leq k-2$. $H_{n-p,n-p-1}^{k,p,l}$ has seven Hamiltonian paths, denoted by $R_1^*,\ldots, R_7^*$, obtained from  Hamiltonian paths $R_1$, $R_3$, and $R_5$ in $H_{n-p,n-p}^{k,p,l}$
by some vertex and edge operations. We present them as follows::
\begin{eqnarray*}
  R_1^*&=& R_1-v_{1t}-u_{1t}v_{1,t-1}+\{u_{21}v_{1,t-1},u_{1t}v_{1,n-k}\}, \\
 R_2^* &=&R_1-v_{1,t+g}-u_{2,g}v_{1,t+g-1}+\{u_{3,s+1}^{(1)}v_{1,t+g-1},u_{2g}v_{1,n-k}\}, \\
  R_3^* &=& R_1-v_{1,n-k},\\
  R_4^* &=& R_3-v_{1,t+g}-u_{2,g}v_{1,t+g-1}+\{u_{11}v_{1,t+g-1},u_{2g}v_{1,n-k}\}, \\
  R_5^* &=& R_3-v_{1,n-k}, \\
  R_6^* &=& R_5-v_{1,t+g+s}-u_{3s}^{(1)}v_{1,t+g+s-1}+\{u_{3,s+1}^{(2)}v_{1,t+g+s-1},u_{3s}^{(1)}v_{1,t+g}\}, \\
  R_7^* &=& R_5-v_{1,t+g+hs}-u_{3s}^{(h)}v_{1,t+g+hs-1}+\{u_{11}v_{1,t+g+hs-1},u_{3s}^{(h)}v_{1,t+g}\}.
\end{eqnarray*}
%
%
%
%
%
%
Hence $H_{n-p,n-p-1}^{k,p,l}$ is Hamilton-biconnected.
Thus $F_{n,n-1}^{k,p,l}$ is  $2p$-Hamilton-biconnected for $ p< l\leq k-2$. Similarly we can prove that
$F_{n,n}^{k,p,l}$ is also $2p$-Hamilton-biconnected for $ p<l= k-1$.

\medskip
We next assume  that $l\leq p$.  Let $r_i\geq 0$ with $\sum_{i=1}^h r_i=p-l$ for $1\leq i\leq h$. Let  $h=k-p$,  $s_0=0$, $s_i=k-l-r_i-1$ for $1\leq i \leq h$, and $t=n-(k-p)(k-l)-l$.  Let  $H_{n-p,m-p}^{k,p,l}(s_1,\ldots,s_{h})$ be a bipartite graph obtained from $F_{n,m}^{k,p,l}$ by deleting all vertices in a balanced set of size $2p$ which consists of vertices with  as large as possible degree (see Fig.~2).
  Let $\mathcal{G}_{n-p,m-p}^{k,p,l}$ be a set of all bipartite graphs $G$ satisfying $G= H_{n-p,m-p}^{k,p,l}(s_1,\ldots,s_{h})$.
Note that every graph of order $m+n-2p$  obtained from $F_{n,m}^{k,p,l}$ by deleting all vertices in a balanced set of size $2p$ contains 
a bipartite graph $G\in \mathcal{G}_{n-p,m-p}^{k,p,l}$ as a subgraph.  It suffices to prove that any bipartite graph  $G\in\mathcal{G}_{n-p,m-p}^{k,p,l}$ is  Hamilton-biconnected.
For any bipartite graph $G\in\mathcal{G}_{n-p,m-p}^{k,p,l}$, without loss of generality, say $G=H_{n-p,m-p}^{k,p,l}(s_1,\ldots,s_{h})$.
 Label the vertices of  $H_{n-p,m-p}^{k,p,l}(s_1,\ldots,s_{h})$ as $u_{11},\ldots,u_{1t}$, $u_{31}^{(1)},\ldots,u_{3,s_1+1}^{(1)},\ldots, u_{31}^{(h)},\ldots,u_{3,s_h+1}^{(h)}$, $v_{11},\ldots,v_{1,m-k}$, $v_{21},\ldots,v_{2h}$ (see Fig.~2).

\medskip
Let $m=n$.  Since $k\geq p+2$ and $0\leq r_i\leq p-l$, we have $h\geq2$ and $s_i\geq1$ for $1\leq i\leq h$.
 Denote
\begin{eqnarray*}
  P_{u_{11}v_{1t}} &=& \bigsqcup_{i=1}^t u_{1i}v_{1i}, \\
  P_{u_{31}^{(h)}v_{2h}} &=& \bigsqcup_{i=1}^{s_h}
 u_{3i}^{(h)}v_{1,t+i+\sum_{j=1}^{h-1}s_j}\bigsqcup u_{3,s_h+1}^{(h)}v_{2h}, \\
  Q_{u_{31}^{(h)}v_{2h}} &=& \bigsqcup_{i=1}^{s_h-1} u_{3i}^{(h)}
 v_{1,t+i+\sum_{j=1}^{h-1}s_j}\bigsqcup u_{3s_h}^{(h)}v_{2h}, \\
  P_{u_{3,s_i+1}^{(i)}v_{1,t+\sum_{j=1}^{i}s_j}} &=&  u_{3,s_i+1}^{(i)}v_{2i}\bigsqcup
\bigg( \bigsqcup_{j=1}^{s_i} u_{3j}^{(i)}v_{1,t+j+\sum_{w=1}^{i-1}s_w}\bigg)~ \mbox{for $1\leq i\leq h$}.
\end{eqnarray*}
$H_{n-p,n-p}^{k,p,l}(s_1,\ldots,s_{h})$ has five kinds of Hamiltonian paths, denoted by $R_1,\ldots, R_5$. We present them as follows:
\begin{eqnarray*}
 R_1&=&P_{u_{11}v_{1t}}\bigsqcup \bigg(\bigsqcup_{i=1}^{h}P_{u_{3,s_i+1}^{(i)}v_{1,t+\sum_{j=1}^{i}s_j}}\bigg),\\
 R_2&=&P_{u_{11}v_{1t}}\bigsqcup \bigg(\bigsqcup_{i=1}^{h-1}P_{u_{3,s_i+1}^{(i)}v_{1,t+\sum_{j=1}^{i}s_j}}\bigg)\bigsqcup P_{u_{31}^{(h)}v_{2h}} ,\\
  R_3&=&\bigsqcup_{i=1}^{h}P_{u_{3,s_i+1}^{(i)}v_{1,t+\sum_{j=1}^{i}s_j}}\bigsqcup P_{u_{11}v_{1t}},\\
 R_4&=& \bigsqcup_{i=1}^{h-1}P_{u_{3,s_i+1}^{(i)}v_{1,t+\sum_{j=1}^{i}s_j}}\bigsqcup P_{u_{11}v_{1t}} \bigsqcup P_{u_{31}^{(h)}v_{2h}},\\
R_5&=&u_{3,s_h+1}^{(h)}v_{1,n-k}\bigsqcup \bigg(\bigsqcup_{i=1}^{h-1}P_{u_{3,s_i+1}^{(i)}v_{1,t+\sum_{j=1}^{i}s_j}}\bigg)\bigsqcup P_{u_{11}v_{1t}} \bigsqcup Q_{u_{31}^{(h)}v_{2h}}.
\end{eqnarray*}
Hence $H_{n-p,n-p}^{k,p,l}(s_1,\ldots,s_{h})$ is Hamilton-biconnected. Thus $F_{n,n}^{k,p,l}$ is $2p$-Hamilton-biconnected for $l\leq p$.

Let $m=n-1$. $H_{n-p,n-p-1}^{k,p,l}(s_1,\ldots,s_{h})$ has four kinds of Hamiltonian paths, denoted by $R_1^*,\ldots, R_4^* $,  obtained from Hamiltonian paths $R_1$ and $R_3$ in  $H_{n-p,n-p}^{k,p,l}(s_1,\ldots,s_{h})$ by some vertex and edge operations. We present them as follows:
\begin{eqnarray*}
R_1^*&=&R_1-v_{1t}-u_{1t}v_{1,t-1}+\{u_{3,s_1+1}^{(1)}v_{1,t-1},u_{1t}v_{1,n-k}\},\\
R_2^*&=&R_1-v_{1,n-k},\\
R_3^*&=&R_3-v_{1,t+s_1}-u_{3,s_1}^{(1)}v_{1,t+s_1-1}+\{u_{3,s_2+1}^{(2)}v_{1,t+s_1-1},u_{3,s_1}^{(1)}v_{1,t}\},\\
R_4^*&=&R_3-v_{1,t+\sum_{i=1}^h s_i}-u_{3,s_h}^{(h)}v_{1,t-1+\sum_{i=1}^h s_i}+\{u_{11}v_{1,t-1+\sum_{i=1}^h s_i},u_{3,s_h}^{(h)}v_{1,t}\}.
 \end{eqnarray*}
Hence $H_{n-p,n-p-1}^{k,p,l}(s_1,\ldots,s_{h})$ is Hamilton-biconnected. Thus  $F_{n,n-1}^{k,p,l}$ is also $2p$-Hamilton-biconnected  for $l\leq p$.
This completes the proof. \QEDB

\medskip

 \noindent{\bf Proof of Lemma~\ref{R2-Lem2.3}.}
(i) Note that every balanced bipartite graph of order $2n-2p$  obtained from $M_{n,n}^{s,t;-}$ by deleting all vertices in a balanced set of size $2p$ contains $M_{s+t+1,s+t+1}^{s,t;-}$ as a subgraph. It suffices to prove that  $M_{s+t+1,s+t+1}^{s,t;-}$ is  Hamilton-biconnected.
Label the vertices of $M_{s+t+1,s+t+1}^{s,t;-}$ as $u_{11},\ldots,u_{1,s-1}$, $u$, $u_{21},\ldots,u_{2,t+1}$, $v_{11},\ldots,v_{1s}$, $v_{21},\ldots,v_{2t}$ (see Fig.~2).
Denote
$$P_{u_{11}v_{1,s-1}}=\bigsqcup_{i=1}^{s-1} u_{1i}v_{1i}, \quad P_{u_{21}v_{2t}}=\bigsqcup_{i=1}^t u_{2i}v_{2i}.$$
 $M_{s+t+1,s+t+1}^{s,t;-}$ has  nine kinds of  Hamiltonian paths, denoted by $R_1,\ldots, R_9$. We present them as follows:
\begin{eqnarray*}
 R_1&=&P_{u_{11}v_{1,s-1}}\bigsqcup uv\bigsqcup P_{u_{21}v_{2t}}\bigsqcup u_{2,t+1}v_{1s},\\
 R_2&=&P_{u_{11}v_{1,s-1}}\bigsqcup P_{u_{21}v_{2t}}\bigsqcup u_{2,t+1}v_{1s}uv,\\
 R_3&=&P_{u_{11}v_{1,s-1}}\bigsqcup uvu_{2,t+1}v_{1s}\bigsqcup P_{u_{21}v_{2t}},\\
 R_4&=&uv\bigsqcup P_{u_{21}v_{2t}}\bigsqcup u_{2,t+1}v_{1s}\bigsqcup P_{u_{11}v_{1,s-1}},\\
 R_5&=&uv_{1,s}\bigsqcup P_{u_{11}v_{1,s-1}}\bigsqcup P_{u_{21}v_{2t}}\bigsqcup u_{2,t+1}v,\\
 R_6&=&uvu_{2,t+1}v_{1s}\bigsqcup P_{u_{11}v_{1,s-1}}\bigsqcup P_{u_{21}v_{2t}},\\
 R_7&=&P_{u_{21}v_{2t}}\bigsqcup u_{2,t+1}vuv_{1s}\bigsqcup P_{u_{11}v_{1,s-1}},\\
 R_8&=&P_{u_{21}v_{2t}}\bigsqcup u_{2,t+1}v_{1s}\bigsqcup P_{u_{11}v_{1,s-1}}\bigsqcup uv,\\
 R_9&=&u_{2,t+1}v_{1s}\bigsqcup P_{u_{11}v_{1,s-1}}\bigsqcup uv\bigsqcup P_{u_{21}v_{2t}}.
  \end{eqnarray*}
Hence $M_{n-p,n-p}^{s,t;-}$ is Hamilton-biconnected. Thus $M_{n,n}^{s,t;-}$ is $2p$-Hamilton-biconnected.

\medskip
(ii) Note that every balanced bipartite graph of order $2n-2p$  obtained from $N_{n,n}^{p,2}$ by deleting all vertices in a balanced set of order $2p$ contains $N_{n-p,n-p}^{0,2}$ as a subgraph. It suffices to prove that  $N_{n-p,n-p}^{0,2}$  is  Hamilton-biconnected.   Let $t=n-p-3$ and label the vertices of $N_{n-p,n-p}^{0,2}$  as $u_{11},\ldots,u_{1t}$, $u_{21},u_{22},u$, $v_{11},\ldots,v_{1t}$, $v_{21},v_{22},v$ (see Fig.~2).
Denote  $$P_{u_{11}v_{1,t-1}}=\bigsqcup_{i=1}^{t-1} u_{1i}v_{1i}, \quad P_{u_{11}v_{1t}}=\bigsqcup_{i=1}^t u_{1i}v_{1i}.$$
$N_{n-p,n-p}^{0,2}$ has nine kinds of Hamiltonian paths. We present them as follows:
\begin{eqnarray*}
 R_1&=&P_{u_{11}v_{1,t-1}}\bigsqcup  u_{2,1}vu_{22}v_{22}uv_{21}u_{1t}v_{1t},\\
  R_2&=&P_{u_{11}v_{1t}}\bigsqcup  u_{2,1}vu_{22}v_{22}uv_{21},\\
 R_3&=&P_{u_{11}v_{1t}}\bigsqcup  u_{2,1}v_{22}uv_{21}u_{22}v,\\
 R_4&=&  u_{2,1}vu_{22}v_{22}uv_{21} \bigsqcup P_{u_{11}v_{1t}},\\
  R_5&=&u_{2,1}vu_{22}v_{1t}\bigsqcup P_{u_{11}v_{1,t-1}} \bigsqcup u_{1t}v_{22}uv_{21},\\
 R_6&=&u_{2,1}v_{22}uv_{21}\bigsqcup P_{u_{11}v_{1t}} \bigsqcup u_{22}v,\\
 R_7&=&uv_{21}u_{22}vu_{21}v_{22}\bigsqcup P_{u_{11}v_{1t}},\\
 R_8&=&uv_{21}u_{22}vu_{21}v_{1t}\bigsqcup P_{u_{11}v_{1,t-1}} \bigsqcup u_{1t}v_{22},\\
 R_9&=&uv_{21}\bigsqcup P_{u_{11}v_{1t}} \bigsqcup u_{21}v_{22}u_{22}v .
 \end{eqnarray*}
Hence $N_{n-p,n-p}^{0,2}$ is Hamilton-biconnected. Thus $N_{n,n}^{p,2}$ is $2p$-Hamilton-biconnected.\QEDB


\bigskip

 \noindent{\bf Proof of Lemma~\ref{R3-Lem2.3}.}
 (i)  Denote $M_{n-p,n-p-1}^{s,t}=(X,Y;E)$ with $|X|=n-p$ and  $|Y|=n-p-1$. Let $x,y\in X$ such that $d(x)=d(y)=n-p-1$.
 Since $s\geq n-t-p-1$, $M_{n-p,n-p-1}^{s,t}$ has no Hamiltonian path between $x$ and $y$. Hence $M_{n-p,n-p-1}^{s,t}$ is not Hamilton-biconnected.
Note that $M_{n-p,n-p-1}^{s,t}$ is one of    graphs   obtained from $M_{n,n-1}^{s,t}$  by deleting all vertices in a balanced set of size $2p$.
It follows from definition that $M_{n-p,n-p-1}^{s,t}$ is not $2p$-Hamilton-biconnected.

\medskip
(ii) Denote $M_{n-p,n-p}^{s,t}=(X,Y;E)$  with $|X|=|Y|=n-p$.   Let $x\in X$ and  $y\in Y$ such that $d(x)=n-p-t$ and  $d(y)=n-p$.
Since $s=n-p-t$,  $M_{n-p,n-p}^{s,t}$  has no Hamiltonian path between $x$ and $y$. Hence  $M_{n-p,n-p}^{s,t}$ is not Hamilton-biconnected.
Note that $M_{n-p,n-p}^{s,t}$ is  one of    graphs  obtained from $M_{n,n}^{s,t}$  by deleting all vertices in  a balanced set of size $2p$.
It follows from definition that $M_{n-p,n-p}^{s,t}$ is not $2p$-Hamilton-biconnected.

\medskip
(iii) Denote $N_{n-p,n-p}^{0,1}=(X,Y;E)$ with $|X|=|Y|=n-p$. Let $x\in X$ and $y\in Y$ such that $d(x)=d(y)=n-p$.
Then $N_{n-p,n-p}^{0,1}$ has no Hamiltonian path between $x$ and $y$.
Hence $N_{n-p,n-p}^{0,1}$ is not Hamilton-biconnected.
Note that $N_{n-p,n-p}^{0,1}$ is one of    graphs  obtained from $N_{n,n}^{p,1}$  by deleting all vertices in a balanced set of size $2p$.
It follows from definition that $N_{n,n}^{p,1}$ is not $2p$-Hamilton-biconnected.\QEDB

\medskip

\noindent {\bf Proof of Lemma~\ref{Lem4.3}.}
 By a similar argument to  Lemma~\ref{Lem3.2}~($i$), $\rho(M_{n,n-1}^{n-k-1,k-p})$, $\rho(M_{n,n-1}^{k-p,n-k-1})$, and $\rho(M_{n,n-1}^{n-k,k-p-1})$
  are the largest roots of $f(x)=0$,  $g(x)=0$, and $h(x)=0$ respectively, where
  $$f(x)=x^4-\Big(n^2-(k-p+1)n+(k+1)(k-p)\Big)x^2+(n-k-1)(n-k+p-1)(k+1)(k-p) ,$$
  $$g(x)=x^4-\Big(n^2-(k-p+1)n+(k-p)(k+1)\Big)x^2+(n-k+1)(n-k+p)k(k-p),$$
  $$h(x)=x^4-\Big(n^2-(k-p)n+k(k-p-1)\Big)x^2+(n-k)(n-k+p)k(k-p-1).$$

  \medskip

(i) Since for all real number $x$,  $$f(x)-g(x)=(n-k-1)(n-2k+p-1)(k-p)>0,$$
 we have $\rho(M_{n,n-1}^{k-p,n-k-1})>\rho(M_{n,n-1}^{n-k-1,k-p})$.

  \medskip

(ii) Since for all real number $x$, $$g(x)-h(x)=(n-2k+p)(x^2+kn-k^2+kp)>0,$$
 we have $\rho(M_{n,n-1}^{n-k,k-p-1})>\rho(M_{n,n-1}^{k-p,n-k-1})$.\QEDB
  \bigskip

\noindent{\bf Proof of Lemma~\ref{Lem4.4}.}
By a similar argument to Lemma~\ref{Lem3.2}~($i$), $q(M_{n,n-1}^{n-k-1,k-p})$, $q(M_{n,n-1}^{k-p,n-k-1})$, and $q(M_{n,n-1}^{n-k,k-p-1})$
  are the largest roots of $f(x)=0$,  $g(x)=0$, and $h(x)=0$, respectively, where $f(x)=xf_1(x)$, $g(x)=xg_1(x)$, and $h(x)=xh_1(x)$,
$$f_1(x)=x^3-\Big(3n+p-1)x^2+(2n^2+(2k+p)n-(2k+1)(k-p+1)\Big)x-(2n-1)(n-k+p-1)(k+1),$$
$$g_1(x)=x^3-\big(3n+p-1)x^2+(2n^2+(2k+p-1)n-k(2k-2p+1)\Big)x-(2n-1)(n-k+p)k, $$
$$h_1(x)=x^3-(3n+p-1)x^2+\Big(2n^2+(2k+p-2)n-(2k-1)(k-p)\Big)x-(2n-1)(n-k+p)k.$$
Since signless Laplacian spectral radius of any nonempty graph is positive, $q(M_{n,n-1}^{n-k-1,k-p})$, $q(M_{n,n-1}^{k-p,n-k-1})$ and $q(M_{n,n-1}^{n-k,k-p-1})$ are the largest roots of $f_1(x)=0$,  $g_1(x)=0$, and $h_1(x)=0$, respectively.

\medskip
(i) 
Since
\begin{eqnarray*}
  f_1(2n-1) &=& (2n-1)(n-k-1)(k-p)>0, \\
  g_1(2n-1) &=& (2n-1)(n-k-1)(k-p)>0,
\end{eqnarray*}
 and for $x\geq 2n-1$,
\begin{eqnarray*}
  f'_1(x) &=& 3x^2-(6n+2p-2)x+2n^2+(2k+p)n-(2k+1)(k-p+1)\\
   &\geq& f'_1(2n-1)\\
   &=& n(n+2k-3p-2)+n^2-(2k+3)(k-p) \\
   &>& (2k-p+2)^2-(2k+3)(k-p) \\
   &\geq&7k+4\\
   &>&0,\\
 g'_1(x) &=& 3x^2-(6n+2p-2)x+2n^2+(2k+p-1)n-k(2k-2p+1)\\
   &\geq& g'_1(2n-1)\\
   &=& n(n+2k-3p-3)+n^2-(k+1)(2k-2p-1) \\
   &>& (2k-p+2)^2-(k+1)(2k-2p-1)  \\
   &\geq&7k+5\\
   &>&0,
\end{eqnarray*}
we have
$$q(M_{n,n-1}^{n-k-1,k-p})<2n-1,\quad q(M_{n,n-1}^{k-p,n-k-1})<2n-1.$$
Together with, for $x< 2n-1$,
$$g_1(x)-f_1(x)=(n-2k+p-1)(2n-1-x)> 0,$$  we have $q(M_{n,n-1}^{n-k-1,k-p})>q(M_{n,n-1}^{k-p,n-k-1})$ for $k\geq p+1$.

(ii) Note that $K_{n,n-k+p-1}$ and $K_{n,n-k+p}$ are proper subgraphs of $M_{n,n-1}^{n-k-1,k-p}$ and $M_{n,n-1}^{n-k,k-p-1}$, respectively.
By Lemma~\ref{Lem2.1},
$$q(M_{n,n-1}^{n-k-1,k-p})>2n-k+p-1, \quad q(M_{n,n-1}^{n-k,k-p-1})>2n-k+p.$$
 Since for $x>2n-k+p-1>n$,
 $$f_1(x)-h_1(x)=(2n-4k+2p-1)x-(2n-1)(n -2k+p-1)>0,$$  we have $q(M_{n,n-1}^{n-k,k-p-1})> q(M_{n,n-1}^{n-k-1,k-p})$.\QEDB
\end{document}